%% file: 01_main.tex
\documentclass{02_iwoloca}
\input{03_preamble_booklet}
\begin{document}

\includepdf{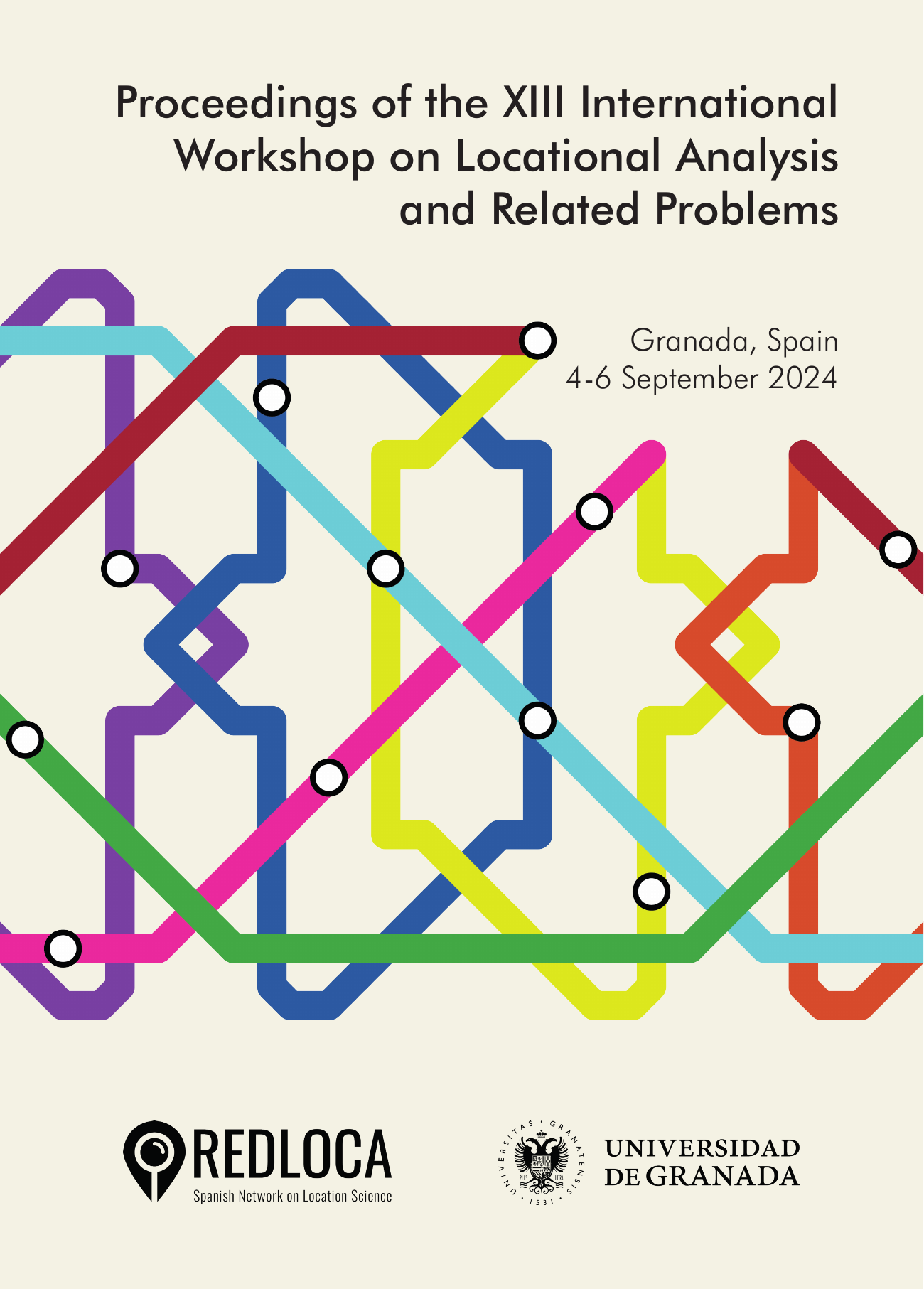}	
\newpage\clearpage\hbox{}\thispagestyle{empty}\newpage
	



\includepdf{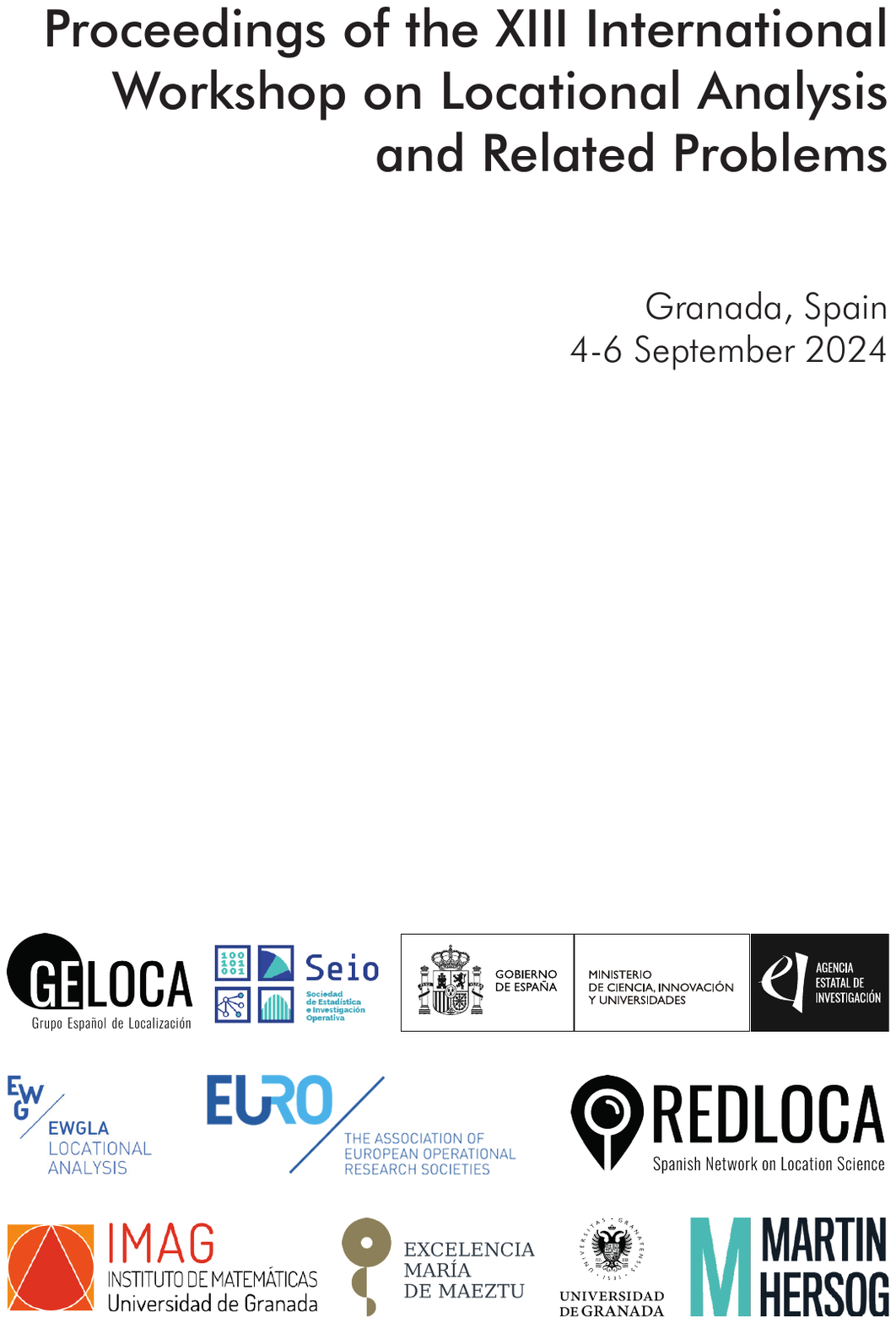}	

\newpage\clearpage\hbox{}\thispagestyle{empty}\newpage
\thispagestyle{empty}
\vspace*{3.5cm}
\begin{flushright}
{\huge PROCEEDINGS OF THE XIII INTERNATIONAL\vspace{-0.2cm}  WORKSHOP ON LOCATIONAL ANALYSIS\\\vspace{0.5cm} AND RELATED PROBLEMS (2024)}
\end{flushright}
\vspace*{4cm}
\begin{flushright}
{
{\small Edited by}\vspace{0.2cm}\\
Marta Baldomero-Naranjo\vspace{0.2cm}\\	
Ricardo G\'azquez\vspace{0.2cm}\\	
Miguel Martínez-Antón\vspace{0.2cm}\\
Luisa I. Martínez-Merino\vspace{0.2cm}\\
Juan M. Muñoz-Ocaña \vspace{0.2cm}\\
Francisco Temprano\vspace{0.2cm}\\
Alberto Torrejón\vspace{0.2cm}\\
Carlos Valverde\vspace{0.2cm}\\
Nicol\'as Zerega
}
\vfill
\texttt{ISBN:  978-84-09-63233-6}
\end{flushright}

\newpage\clearpage\hbox{}\thispagestyle{empty}\newpage

\tableofcontents


\input{04_about}

\chapter{Program}

I: Invited Talk, C: Contributed Talk
\input{05_timetable.tex}

\cleardoublepage

\chapter{Invited Speakers}
\thispagestyle{empty}

\newpage\clearpage\thispagestyle{empty}\hbox{}\newpage
\input{abstracts/plenaries/1.marin} 

\newpage\clearpage\thispagestyle{empty}\hbox{}\newpage
\input{abstracts/plenaries/2.luedtke-iwoloca} 

\newpage
\input{abstracts/plenaries/3.salman} 

\chapter{Abstracts}
\thispagestyle{empty}


\newpage\clearpage\thispagestyle{empty}\hbox{}\newpage

\input{abstracts/tex/IWOLOCA24_Abstract_9444} 
\newpage

\input{abstracts/tex/IWOLOCA24_Abstract_1691} 
\newpage\clearpage\thispagestyle{empty}\hbox{}\newpage

\input{abstracts/tex/IWOLOCA24_Abstract_8661} 
\newpage

\input{abstracts/tex/IWOLOCA24_Abstract_0233} 
\newpage

\input{abstracts/tex/IWOLOCA24_Abstract_6905} 
\newpage

\input{abstracts/tex/IWOLOCA24_Abstract_5013} 
\newpage\clearpage\thispagestyle{empty}\hbox{}\newpage

\input{abstracts/tex/IWOLOCA24_Abstract_4203} 
\newpage

\input{abstracts/tex/IWOLOCA24_Abstract_9068} 
\newpage

\input{abstracts/tex/IWOLOCA24_Abstract_7341} 
\newpage

\input{abstracts/tex/IWOLOCA24_Abstract_5262} 
\newpage

\input{abstracts/tex/IWOLOCA24_Abstract_2207} 
\newpage

\input{abstracts/tex/IWOLOCA24_Abstract_9914} 
\newpage

\input{abstracts/tex/IWOLOCA24_Abstract_5122} 
\newpage

\input{abstracts/tex/IWOLOCA24_Abstract_2696} 
\newpage

\input{abstracts/tex/IWOLOCA24_Abstract_2643} 
\newpage

\input{abstracts/tex/IWOLOCA24_Abstract_2489} 
\newpage

\input{abstracts/tex/IWOLOCA24_Abstract_7279} 
\newpage\clearpage\thispagestyle{empty}\hbox{}\newpage

\input{abstracts/tex/IWOLOCA24_Abstract_7300} 
\newpage

\input{abstracts/tex/IWOLOCA24_Abstract_9592} 
\newpage

\input{abstracts/tex/IWOLOCA24_Abstract_5088} 
\newpage

\input{abstracts/tex/IWOLOCA24_Abstract_8729} 
\newpage\clearpage\thispagestyle{empty}\hbox{}\newpage

\input{abstracts/tex/IWOLOCA24_Abstract_2730} 
\newpage

\input{abstracts/tex/IWOLOCA24_Abstract_1386} 
\newpage\clearpage\thispagestyle{empty}\hbox{}\newpage

\input{abstracts/tex/IWOLOCA24_Abstract_3749} 
\newpage\clearpage\thispagestyle{empty}\hbox{}\newpage

\input{abstracts/tex/IWOLOCA24_Abstract_7463} 
\newpage\clearpage\thispagestyle{empty}\hbox{}\newpage

\input{abstracts/tex/IWOLOCA24_Abstract_4806} 


\cleardoublepage 
\printindex[author]
 
\chapter{Useful Information}

\input{06_useful_info}

\chapter{Partner Institutions and Sponsors}

\input{07_logos.tex}

\newpage


\includepdf{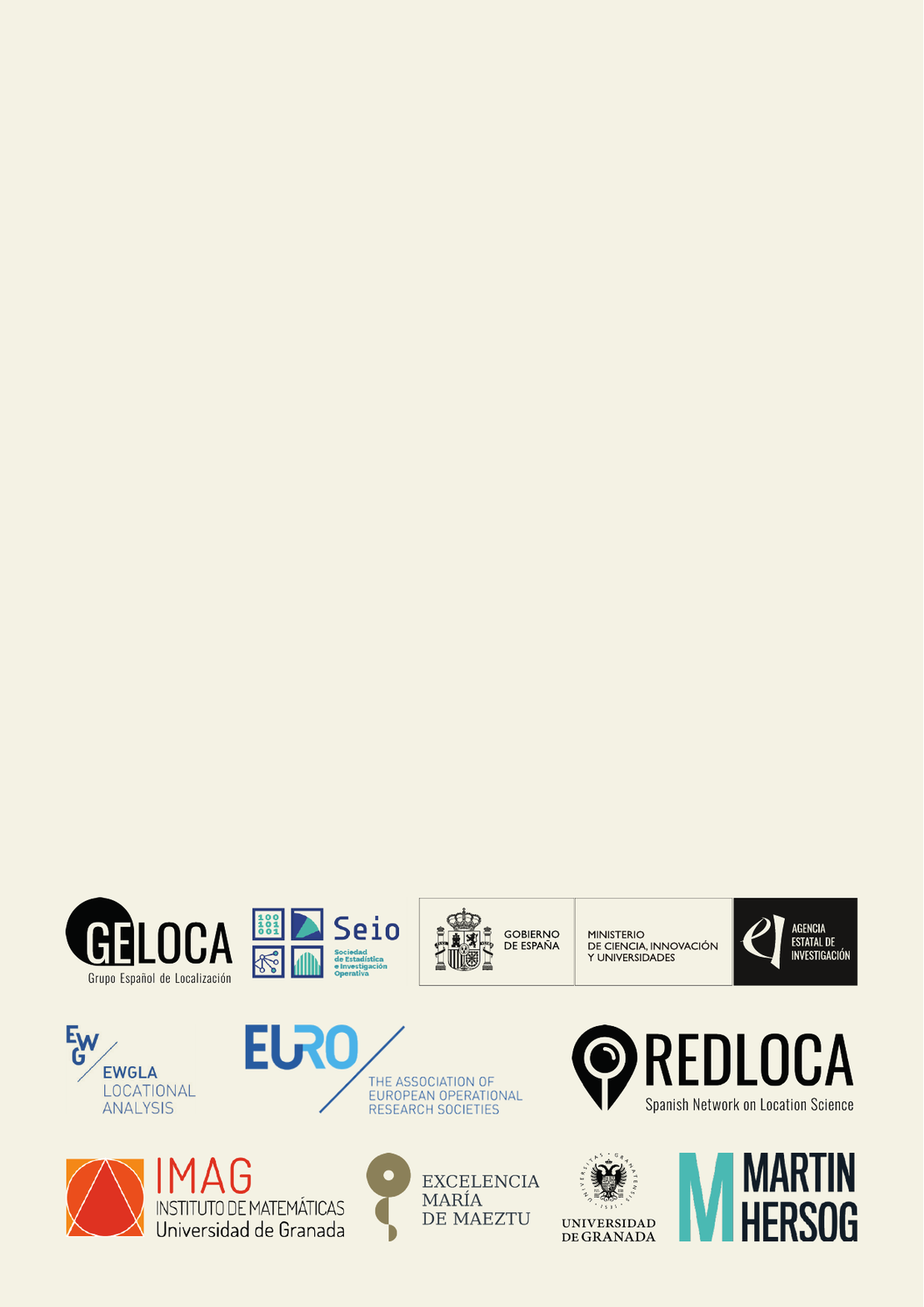}	

\end{document}

%% file: 03_preamble_booklet.tex
\definecolor{myorange}{RGB}{255,117,40}
\definecolor{mygray}{RGB}{164, 168, 172}
\definecolor{mywhite}{RGB}{235, 238, 231}
\definecolor{myblue}{RGB}{29, 87, 152}

\newenvironment{abstract_online}[4] 
{\filbreak 
	
	{\Large \bfseries #1} \\
	
	{\bfseries \itshape #2} \\
    {#3}    
    \vspace*{-0.25cm}
    \hrule
    \vspace*{0.25cm}
	\textcolor{mygray}{\textbf{Keywords:} #4}
	\vspace*{0.5cm}
}
{}

\usepackage{etoolbox}
\usepackage[scale=1,angle=0,opacity=1]{background}
\backgroundsetup{contents={}}

\newdimen\mybarpadding
\RedeclareSectionCommand[%
    ,afterskip=4em plus 1pt minus 1pt%
    ,beforeskip=-1pt
    ,level=0%
    ,toclevel=0%
]{chapter}%

\setkomafont{chapter}{\normalfont\normalsize\bfseries\Huge} 

\newcommand*{\mynumberedtest}[1]{
  \if\relax\detokenize{#1}\relax%
  \else%
    #1%
    
  \fi}

\makeatother

\newcommand{\tbg}{gray} 
\newcommand{\tfg}{white}
\newcommand{\tbc}{gray!25}

\newcommand{\IScolor}{myblue!65} 
\newcommand{\CTcolor}{white} 
\newcommand{\tablebreak}[2]{
	\rowcolor{\tbc} #1 &  \multicolumn{4}{c|}{\bfseries #2} \\ \hline }
\newcommand{\eventtype}[2]{
	#1& \multicolumn{4}{c|}{\cellcolor{\tbg}\color{\tfg}\bfseries #2} \\ \hline }

\newcolumntype{L}[1]{%
	>{\raggedright\let\newline\\\arraybackslash\hspace{0pt}}m{#1}}
\newcolumntype{C}[1]{%
	>{\centering\let\newline\\\arraybackslash\hspace{0pt}}m{#1}}
\newcolumntype{R}[1]{%
	>{\raggedleft\let\newline\\\arraybackslash\hspace{0pt}}m{#1}}

\newcommand{\IS}[5]{
	#1 &\cellcolor{\IScolor}I&{\bfseries#2}\newline #4&&#5 \\ \hline}
\newcommand{\CT}[5]{%
	#1 &\cellcolor{\CTcolor}C&{\bfseries#2}\newline #4&&#5 \\ \hline}

\usepackage[overwrite=true,mode=gm]{getmap}

%% file: 04_about.tex
\chapter{Welcome}
\thispagestyle{plain}
Dear Locators,

Welcome to the XIII International Workshop on Locational Analysis and Related Problems (IWOLOCA 2024) in Granada, a meeting organized by the Spanish Network on Locational Analysis (REDLOCA) and the Location Group GELOCA (Spanish Society of Statistics and Operations Research--SEIO). Our workshops aim to bring together researchers and
practitioners working in the wide area of Location Science through the proposal of
models and methods to solve theoretical and practical problems. These meetings
traditionally combine a strong scientific value with a very enjoyable and friendly
come together as an enthusiastic and energetic group, where old friends meet again
and new members are heartily welcomed.

The success of these meetings over the years proves the continuing research community interested in locational analysis and related problems. This 13th edition of the workshop features 26 
contributions covering a broad range of topics, from discrete location, over hub location and routing, to
network design and several interesting applications, including those in cross-docking designs. During the three days of the meeting the
works will be presented in 7 sessions, with around 40 participants. The workshop will feature the participation of three highly recognized plenary speakers who accepted our invitation to come to Granada and share their insights with us. It is a pleasure to have them in Granada, thank you! Prof. Mar\'in (Universidad de Murcia) will present his work on an interesting application of location problems to determining original DNA sequences from different fragments. Prof. Luedtke (UW-Madison) will bring models and solution methodologies for solving a location problem that arises when designing power generation plants based on the impact of extreme weather events. Eventually, Prof. Salman (Ko\c{c} University) will share different studies on disaster management and locational decisions with us.

We would like to thank all the members of the Organizing Committee for their valuable help in the organization, and we are very grateful to the colleagues of the Scientific Committee. They acknowledge the editors of these proceedings, for the hard work committed to collecting all the information about the workshop and giving it this fantastic shape. We would also like to thank our graphical designer, Mar\'ia Mart\'in Hersog, for her help in designing all the amazing drawings identifying the workshop in this edition. Finally, we would like to thank the support of the Spanish Agency of Research (AEI) through grants RED2022-134149-T funded
by MICIU/AEI /10.13039/501100011033 and IMAG-Maria de Maeztu grant CEX2020-
001105-M /AEI /10.13039/501100011033, Universidad de Granada, GELOCA (SEIO), and Institute of Mathematics (IMAG) as well.

Our workshop represents a great opportunity to build new and enrich existing relationships and to share experiences among locators. We wish you all a successful and fruitful meeting, with new ideas and
collaborations, and a pleasant stay in Granada.

\vspace*{1cm}

V\'ictor Blanco\\
IWOLOCA 2024 Organizing Committee, chair

\newpage

\chapter{IWOLOCA}

The XIII International Workshop on Locational Analysis and Related Problems will take place during September 4--6, 2024 in Granada (Spain). It is organized by the Spanish Location Network, REDLOCA, and the Location Group, GELOCA, from the Spanish Society of Statistics and Operations Research (SEIO). The Spanish Location Network is a group of 100+ researchers from several Spanish universities organized into 7 thematic groups. The Network has been funded by the Spanish Government since 2003. 

The topics of interest are location analysis and related problems. This includes location models, networks, transportation, logistics, exact and heuristic solution methods, and computational geometry, among many others.

\section*{Previous meetings}
One of the main activities of the Network is a yearly meeting aimed at promoting communication among its members and between them and other researchers and contributing to the development of the location field and related problems. The last meetings have taken place in Edinburgh (September 7--8, 2023), Elche (January 31--February 1, 2022), Sevilla (January 23--24, 2020), C\'adiz (January 20--February 1, 2019), Segovia (September 27--29, 2017), Málaga (September 14--16, 2016),  Barcelona (November 25--28, 2015), Sevilla  (October 1--3, 2014),  Torremolinos (M\'alaga, June 19--21, 2013),  Granada (May 10--12, 2012),   Las Palmas de Gran Canaria  (February 2--5, 2011) and Sevilla (February 1--3, 2010).

\chapter{Committees}


\begin{center}
\begin{tabular}{ll}
\multicolumn{2}{l}{\bf \large Scientific Committee}\\
\multicolumn{2}{l}{}\\
Mari Albareda Sambola   & Universidad Politécnica de Cataluña \\
Víctor Blanco  & Universidad de Granada \\
David Canca  & Universidad de Sevilla \\
Sergio García Quiles  & University of Edinburgh \\
Jörg Kalcsics  & University of Edinburgh \\
Mercedes Landete  & Universidad Miguel Hernández de Elche \\
Teresa Ortuño  & Universidad Complutense de Madrid \\
Blas Pelegrín  & Universidad de Murcia \\
Justo Puerto Albandoz  & Universidad de Sevilla \\
Antonio M. Rodríguez-Chía  & Universidad de Cádiz \\
Juan José Salazar  & Universidad de La Laguna \\
\multicolumn{2}{l}{}\\
\multicolumn{2}{l}{}\\
\multicolumn{2}{l}{\bf \large Organizing Committee}\\
\multicolumn{2}{l}{}\\
Víctor Blanco & Universidad de Granada \\
Ricardo Gázquez & Universidad Carlos III de Madrid\\
Gabriel González & Universidad de Granada\\
Miguel Martínez-Antón & Universidad de Granada\\
Rom\'an Salmer\'on & Universidad de Granada
\end{tabular}
\end{center}

%% file: 05_timetable.tex

\section{Wednesday, September 4, 2024}
\begin{center}
	\filbreak
\begin{longtable}{|C{0.15\linewidth}| C{0.04\linewidth}|  C{0.3\linewidth} C{0.0\linewidth} C{0.4\linewidth}|}\hline

    \tablebreak{4:20--4.30}{Registration}
    \tablebreak{4:30--4:45}{Opening session}
    
    \eventtype{4:45--5:45}{Plenary talk I}
        \IS{}{Alfredo Marín}{}{U. Murcia, Spain}{A new application of discrete location: DNA sequence assembly}
    
    \tablebreak{5:45--6:00}{Coffee}
    
    \eventtype{6:00--7:40}{Session 1: Applications}
        \CT{6:00--6:25}{Alberto Japón}{}{U. Sevilla, Spain}{Classification forests via mathematical programming: case study
        in obesity detection}
        \CT{6:25--6:50}{Ramón Piedra-de-la-Cuadra}{}{U. Sevilla, Spain}{Optimization approaches to scheduling working shifts for train dispatchers}
        \CT{6:50--7:15}{Eduardo Pipicelli}{}{U. Naples Federico II, Italy}{A mathematical model for designing a postal banking network}
        \CT{7:15--7:40}{David Canca (chair)}{}{U. Sevilla, Spain}{A bi-level metaheuristic for the design of hierarchical transit
        networks in a multimodal context}

    \tablebreak{8:30}{Welcome Cocktail: \textbf{BHeaven Granada}}
\end{longtable}
\end{center}

\newpage\clearpage\thispagestyle{empty}\hbox{}\newpage

\section{Thursday, September 5, 2024}
\begin{center}
\filbreak
\begin{longtable}{|C{0.15\linewidth}| C{0.04\linewidth}|  C{0.3\linewidth} C{0.0\linewidth} C{0.4\linewidth}|}\hline	

    \eventtype{9:00--10:00}{Plenary talk II}
        \IS{}{Jim Luedtke
        }{}{U. Wisconsin-Madison, USA}{Generator Expansion to Improve Power Grid Resiliency and Efficiency: A Case Study in Location Analysis Under Uncertainty}

    \eventtype{10:00--11:15}{Session 2: Covering}
        \CT{10:00--10:25}{Marta Baldomero 
        Naranjo}{}{U. Cádiz, Spain}{The maximal covering location problem with edge downgrades}
        \CT{10:25--10:50}{Ricardo Gázquez
        }{}{U. Carlos III Madrid, Spain}{An upgrading approach for the multi-type Maximal Covering Location Problem}
        \CT{10:50--11:15}{Antonio Rodríguez Chía (chair)}{}{U. Cádiz, Spain}{New results in the covering tour problem with edge upgrades}
    
    \tablebreak{11:15--11:35}{Coffee}
    \tablebreak{11:35--11:50}{Groups meetings}
    
    \eventtype{11:50--1:30}{Session 3: Discrete Location}
        \CT{11:50--12:15}{Concepción Domínguez}{}{U. Murcia, Spain}{Stable Solutions for the Capacitated Simple Plant Location Problem with Order}
        \CT{12:15--12:40}{Sophia Wrede}{}{RWT Aachen U., Germany,}{Cover-based inequalities for the single-source capacitated facility location problem with customer preferences}
        \CT{12:40--1:05}{Juan José Salazar González}{}{U. La Laguna, Spain}{Solving a multi-source capacitated facility location problem with
        a fairness objective}
        \CT{1:05--1:30}{José Fernando Camacho Vallejo (chair)}{}{Tec Monterrey, Mexico}{A Reformulation for the Medianoid Problem with Multipurpose

        Trips}
\end{longtable}
\end{center}

\newpage
\begin{center}
	\filbreak
\begin{longtable}{|C{0.15\linewidth}| C{0.04\linewidth}|  C{0.3\linewidth} C{0.0\linewidth} C{0.4\linewidth}|}\hline	
    
    \tablebreak{1:30--3:15}{Lunch}
    
    \eventtype{3:15--4:30}{Session 4: Network I}
        \CT{3:15--3:40}{Francisco Temprano}{}{U. Sevilla, Spain}{New advances in hypergraph structure analysis}
        \CT{3:40--4:05}{Gabriel González
        }{}{U. Granada, Spain}{Expanding Reaction Networks through Autocatalytic Subnetworks
        that Maximize Growth Factor}
        \CT{4:05--4:30}{Carlos Valverde (chair)}{}{U. Sevilla, Spain}{The Hampered K-Median Problem with Neighbourhoods}
    
    \eventtype{4:30--5:45}{ Session 5:  Hub Location and Routing}
        \CT{4:30--4:55}{Juan Manuel Muñoz-Ocaña}{}{U. Cádiz, Spain}{Formulations and resolution procedures for upgrading hub networks}
        \CT{4:55--5:20}{Inmaculada Espejo}{}{U. Cádiz, Spain}{A further research on Stochastic Single-Allocation Hub Location Problem}
        \CT{5:20--5:45}{Mari Albareda (chair)}{}{U. Polit\'ecnica Catalunya, Spain}{A single vehicle location routing problem with pickup and delivery}
    
    \tablebreak{5:45}{Group Picture \& Guided Tour}
\end{longtable}
\end{center}

\newpage
\section{Friday, September 6, 2024}
\vspace*{-0.5cm}
\begin{center}
\begin{longtable}{|C{0.15\linewidth}| C{0.04\linewidth}|  C{0.3\linewidth} C{0.0\linewidth} C{0.4\linewidth}|}\hline	

    \eventtype{9:00--10:00}{Plenary talk III}
        \IS{}{Sibel Salman}{}{Ko\c{c} U., Turkey}{Location, Allocation, Routing and Network Design in Humanitarian Operations}
    
    \eventtype{10:00--11:15}{Session 6: Routing}
        \CT{10:00-10:25}{Inigo Martin Melero}{}{U. Miguel Hern\'andez Elche, Spain}{The Effect of Budget Limiting on the Linear Ordering Problem}
        \CT{10:25--10:50}{Aitor López-Sánchez}{}{U. Rey Juan Carlos, Spain}{Synchronization routing for agricultural vehicles and implements}
        \CT{10:50--11:15}{Paula Segura (chair)}{}{U. Valencia, Spain}{Exact approaches for the Chinese Postman Problem with load-dependent costs}
    
    \tablebreak{11:15--11:35}{Coffee}
    \tablebreak{11:35--12:15}{Groups meetings}
    
    \eventtype{12:15--1:30}{Session 7: Combinatorial Optimization}
        \CT{12:15--12:40}{Alberto Torrejón}{}{U. Sevilla, Spain}{Modeling ordered interactions in Location Problems}
        \CT{12:40--1:05}{Miguel Martínez-Antón}{}{U. Granada, Spain}{Finding a smallest mediated set}
        \CT{1:05--1:30}{Martine Labbé (chair)}{}{U. Libre de Bruxelles, Belgium}{Novel Valid Inequalities for Chance-Constrained Problems with Finite Support}		
\end{longtable}
\end{center}

\newpage
\begin{center}
\filbreak
\begin{longtable}{|C{0.15\linewidth}| C{0.04\linewidth}|  C{0.3\linewidth} C{0.0\linewidth} C{0.4\linewidth}|}\hline	
    
    \tablebreak{1:30--3:15}{Lunch}
    
    \eventtype{3:15--4:30}{Session 8: Cross-docking door design}
        \CT{3:15--3:40}{Aitziber Unzueta}{}{U. Pais Vasco, Spain}{Cross-docking platforms design and management under uncertainty}
        \CT{3:40-4:05}{Laureano Escudero}{}{U. Rey Juan Carlos, Spain}{Cross-docking platforms design and distributionally robust optimization}
        \CT{4:05-4:30}{Maria Araceli Garin (chair)}{}{U. Pais Vasco, Spain}{Cross-docking platforms design and  mixed binary quadratic model  for distributionally robust optimization}
    
    \tablebreak{4:30--4:45}{Closing Session}
    \tablebreak{4:45--5:15}{Coffe Break}
    \tablebreak{5:15--5:45}{REDLOCA Meeting: Spanish Network on Locational Analysis}
    \multicolumn{2}{c}{}\\\hline
    \tablebreak{9:00}{Gala Dinner: \textbf{Carmen de la Victoria}}
	\end{longtable}
\end{center}

%% file: abstracts/plenaries/1.marin.tex

\index[author]{Mar\'in, Alfredo \\\hspace{0.5cm}
	 Universidad de Murcia, Spain, \email{amarin@um.es}}


\begin{abstract_online}
{
A new application of discrete location: DNA sequence assembly
}
{
Alfredo Marín
}
{
Universidad de Murcia,  \email{amarin@um.es}. \\
}
{
Facility location, Discrete Location
}
\end{abstract_online}


Given fragments of several copies of a unique DNA sequence, "sequence assembly" is the problem of determining the original sequence from these fragments. Integer Programming formulations coming from the vehicle routing field have been used sometimes to solve the problem. In this talk we try to re-interpret the problem as a kind of discrete location problem and design formulations based on this new interpretation.


%% file: abstracts/plenaries/2.luedtke-iwoloca.tex


\index[author]{Luedtke, James\\\hspace{0.5cm}
	 University of Wisconsin-Madison, USA, \email{jim.luedtke@wisc.edu}}
\index[author]{Rossmann, Ramsey \\\hspace{0.5cm}
	University of Wisconsin-Madison, USA, \email{rossmann2@wisc.edu}}
\index[author]{Anitescu, Mihai \\\hspace{0.5cm}
	 Argonne National Laboratory, USA, \email{anitescu@mcs.anl.gov}}
\index[author]{Bessac, Julie\\\hspace{0.5cm}
	 Argonne National Laboratory, USA, \email{jbessac@anl.gov}}
\index[author]{Ferris, Michael\\\hspace{0.5cm}
	University of Wisconsin-Madison, USA, \email{ferris@cs.wisc.edu}}
\index[author]{Krock, Mitchell\\\hspace{0.5cm}
	 Argonne National Laboratory, USA}
\index[author]{Roald, Line\\\hspace{0.5cm}
	University of Wisconsin-Madison, USA, \email{roald@wisc.edu}}

\begin{abstract_online}
{
Generator Expansion to Improve Power Grid Resiliency and Efficiency: A Case Study in Location Analysis Under Uncertainty
}
{
James Luedtke\,\textsuperscript{a,*}, Ramsey Rossmann\,\textsuperscript{a}, Mihai
Anitescu\,\textsuperscript{b}, Julie Bessac\,\textsuperscript{b}, Mitchell Krock\,\textsuperscript{b}, Line
Roald\,\textsuperscript{a}
}
{
\textsuperscript{a }University of Wiconsin-Madison,  \email{jim.luedtke@wisc.edu}\\
\textsuperscript{b }Argonne National Laboratory,\\ 
\textsuperscript{* }Presenting author. \\
}
{
Facility Location, Power Grid, Resiliency, Stochastic Programming
}
\end{abstract_online}


\section*{Introduction}

We consider the problem of choosing power generator locations in a power grid to yield a system that is both efficient on average and resilient to extreme weather events. This problem is challenging due to the combination of discrete decisions and the need to consider performance of the system under rare high-impact outcomes. We use this problem as a case study to illustrate modeling principles and solution techniques that may be useful more broadly in location analysis problems under uncertainty. Specifically, we introduce a stochastic integer programming modeling framework, advocate for a bi-objective modeling approach for balancing efficiency and resiliency, introduce a conditional sampling method for addressing the challenge of low-probability high-impact events, discuss methods for solving the resulting approximation models, and illustrate the importance of modeling spatial dependence of the uncertain parameters in the system.

\section*{Sample Results}

\begin{figure}[H]
\begin{center}
\includegraphics[scale=0.2]{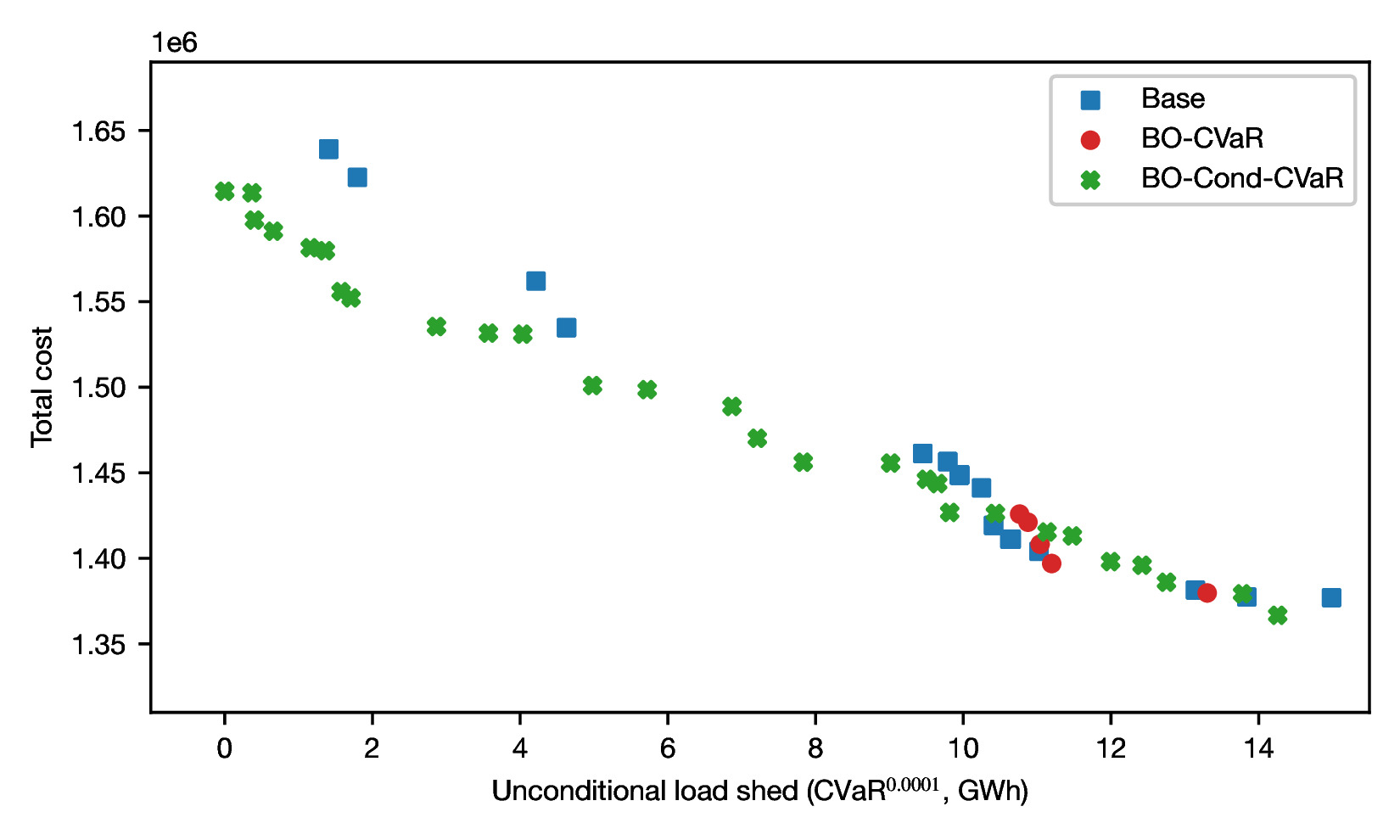}
\end{center}
\caption{Comparison of efficient frontiers obtained from traditional stochastic programming model, bi-objective model,
and bi-objective model with conditional sampling.}\label{fig:ef}
\end{figure}

Figure \ref{fig:ef} displays the results of our proposed model and compares it to alternatives. This figure displays the
two relevant metrics of interest, {\it average cost} and the {\it average load shed} in the 0.01\% highest load-shed hours, which is
a measure of the resiliency of the system to extreme events.  Varying model parametrs allows each method to obtain
solutions that trade-off these two objectives. Model `base' is a traditional two-stage stohastic programming model that
minimizes expected cost, where the cost includes a penalty on load shed. Model `BO-CVaR' is a
bi-objective model that explicitly considers the two objectives of minimum expected cost and minimum load shed in
extreme scenarios, and is solved by standard sample average approximation (SAA). This method fails to find solutions
of low risk due to the challenge in estimating the risk measure in extreme events when using SAA. Model `BO-CVaR-Cond'
is the same model as `BO-CVaR' but adjusts the sampling to generate scenarios from a distribution conditional on the
temperatures being in either the 1\% extreme lowest or highest scenarios. This method is much more capable of generating
solutions with low risk of load shedding. 

\begin{figure}[H]
	\centering
	\begin{subfigure}[b]{.48\linewidth}
		\centering
		\includegraphics[scale=0.28]{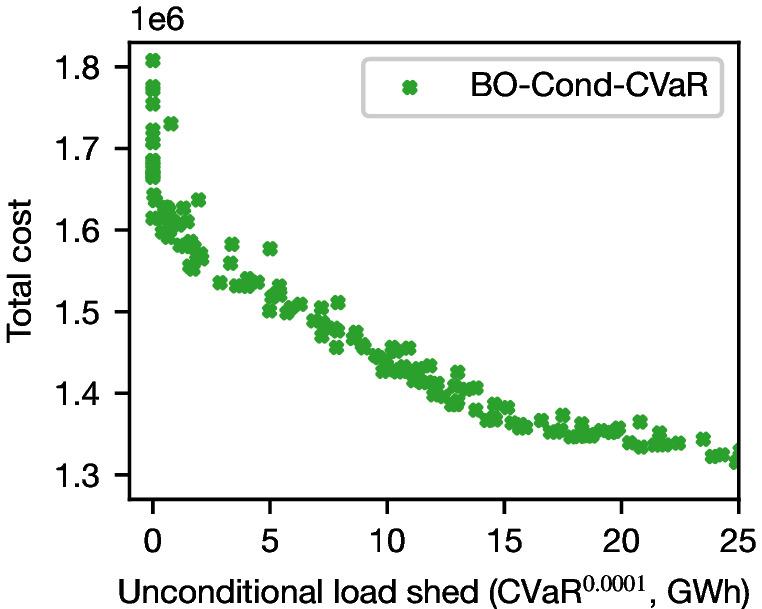}
		\caption{Dependent model.}
	\end{subfigure}~
	\begin{subfigure}[b]{.48\linewidth}
		\centering 
		\includegraphics[scale=0.28]{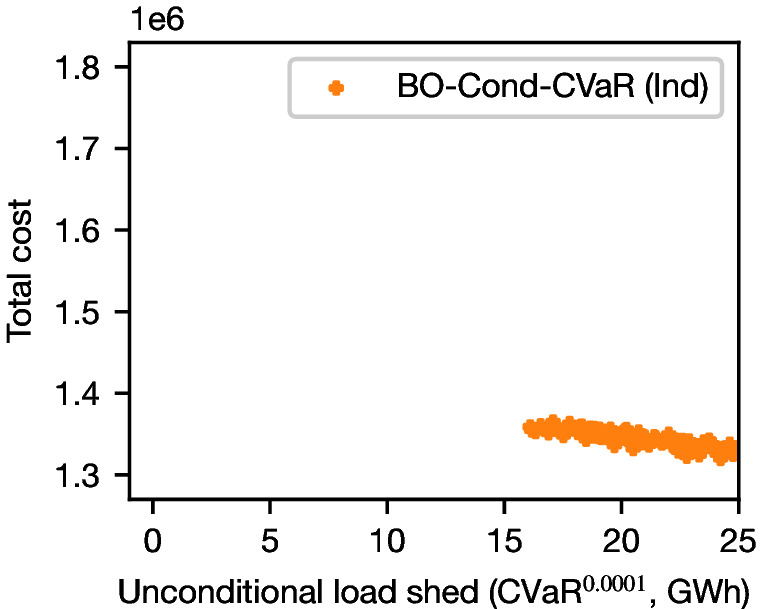}
		\caption{Independent model.}
	\end{subfigure}
	\caption{Comparison of efficient frontiers obtained from model that considers spatial dependence and model that
	assumes spatial independence.}
	\label{fig:dep}
\end{figure}

Figure \ref{fig:dep} illustrates the importance of modeling spatial dependence in our models -- without doing so, the
model is unable to generate solutions with low risk of load shedding.

%% file: abstracts/plenaries/3.salman.tex

\index[author]{Salman, Sibel \\\hspace{0.5cm}
	 Koç University, Turkey, \email{ssalman@ku.edu.tr}}


\begin{abstract_online}
{
Location, Allocation, Routing and Network Design in \\ Humanitarian Operations}
{
Sibel Salman
}
{
Koç University,  \email{ssalman@ku.edu.tr}.  \\
}
{
Humanitarian Operations, Facility Location, Network Design, Vehicle Routing
}
\end{abstract_online}

\section*{Introduction}

Due to the increasing number of people affected by disasters worldwide and the increase in weather-related disasters, developing mathematical models and solution methodology has become increasingly important to reduce risk and ensure efficient response. We will present two studies on disaster preparedness and response: i) relief aid distribution after a natural disaster and ii) providing vaccination to remote areas during a pandemic.

\section*{Designing an efficient and equitable relief aid distribution network with sorting and recycling facilities}
In post-disaster response, relief items are delivered to meet immediate needs and alleviate suffering. Alongside contributions from organizations, individuals often send in-kind donations, which can be unsuitable and complicate management, leading to resource wastage and delays. Establishing sorting facilities can mitigate these issues by ensuring only relevant items reach distribution points, with others stored or recycled. We focus on designing an efficient and equitable relief aid distribution network and investigate the impact of relief aid sorting and reverse supply chain strategies by optimizing a network configuration with and without sorting facilities. The nodes consist of donation centers as supply points, sorting centers as intermediate facilities, points of distribution, relief agency warehouses, and recovery centers as destination points. Road and air shipment options exist between these locations, with helicopters offering faster delivery but at higher costs.

In the case with sorting centers and recycling, donated items are sent to sorting centers, where they are categorized into three types: Type I (urgent items), Type II (useful for future disasters), and Type III (broken or unusable items). Type I items are quickly dispatched to local distribution points, Type II items are stored in relief agency warehouses, and Type III items are sent to recovery centers for processing and revenue generation. We use a three-stage stochastic programming model to optimize the location and shipment decisions. First, sorting center locations are determined pre-disaster under demand and supply uncertainty. Post-disaster, realized demand and travel times dictate shipment quantities to sorting centers. Once sorted, shipment quantities, destinations for each item type, and transportation modes are determined. The goal is to minimize total expected costs, including setup, workforce, shipment, and unsatisfied demand costs, offset by revenue from Type III items. In the unsorted donations case, donations with unknown mixes of Type I, II, and III items are sent directly from donation centers to PODs. At the PODs, items are categorized upon arrival: Type I items go to victims, while Types II and III are discarded. Thus, only a portion of shipped items meets demand. The total supply chain cost includes shipment costs and the penalty for unsatisfied demand. We use a single-stage stochastic programming model to determine shipment quantities from donation centers to PODs. We applied the two models to the anticipated Istanbul earthquake scenario and compared sorted and unsorted donation cases based on total cost and the latest shipment time.

\section*{Routing and Scheduling of Mobile Vaccination Services in a Pandemic}

During a pandemic, the widespread distribution of vaccines is crucial to avoid an explosive increase in infections. However, some people have difficulty accessing vaccination services at fixed locations, such as those living in rural areas, disabled or elderly people who are immobile, and refugees, for which it is observed that vaccination rates are lower. Fixed vaccination services are extended with mobile facilities that get close to the locations of such people. We optimize the routes and schedules of the mobile facilities over a medium-term planning horizon to minimize the number of people that cannot be reached due to capacity restrictions, the total lateness in vaccination, the total distance traveled by the mobile facilities, and inequity over the locations. We formulate a multi-period selective routing model and incorporate several alternative approaches for representing equity. We develop a heuristic solution approach to solve the case of COVID-19 vaccination of Syrian refugee groups, in addition to Turkish citizens living in different neighborhoods of the city of Gaziantep in Turkey. We develop a three-stage matheuristic and a metaheuristic to solve this large-scale instance and derive insights on service quality via various analyses.



%% file: abstracts/tex/IWOLOCA24_Abstract_9444.tex


\index[author]{Albareda, Mar\'ia \\\hspace{0.5cm}
	 Universitat Polit\'ecnica de Catalunya, Spain, \email{maria.albareda@upc.edu}}
\index[author]{Blanco, V\'ictor \\\hspace{0.5cm}
	 Universidad de Granada, Spain, \email{vblanco@ugr.es}}
\index[author]{Hinojosa, Yolanda \\\hspace{0.5cm}
	 Universidad de Sevilla, Spain, \email{yhinojos@us.es}}


\title{A single vehicle location routing problem with pickup and delivery}


\author{
 Mari Albareda\,\textsuperscript{a} , V\'ictor Blanco\,\textsuperscript{b}, Yolanda Hinojosa\,\textsuperscript{c}}

\address{
      \textsuperscript{a }Department d’Estad\'istica i Investigaci\'o Operativa. Universitat Polit\'ecnica de Catalunya, Spain, \email{maria.albareda@upc.edu}\\
      \textsuperscript{b }Institute of Mathematics (IMAG). Dpt. Quant. Methods for Economics \& Business. Universidad de  Granada, Spain, \email{vblanco@ugr.es}\\
      \textsuperscript{c }Institute of Mathematics (IMUS). Dpt. Applied Economics I. Universidad de Sevilla, Spain, \email{yhinojos@us.es}\\
}

\keywords{Facility Location, Network Design, IWOLOCA}

\begin{abstract_online}
{A single vehicle location routing problem with \newline pickup and delivery}
{
 Mari Albareda\,\textsuperscript{a,*}, V\'ictor Blanco\,\textsuperscript{b}, Yolanda Hinojosa\,\textsuperscript{c}
}
{
\textsuperscript{a }Department d’Estad\'istica i Investigaci\'o Operativa. Universitat Polit\'ecnica de Catalunya, Spain, \email{maria.albareda@upc.edu}\\
\textsuperscript{b }Institute of Mathematics (IMAG). Dpt. Quant. Methods for Economics \& Business. Universidad de  Granada, Spain, \email{vblanco@ugr.es}\\
\textsuperscript{c }Institute of Mathematics (IMUS). Dpt. Applied Economics I. Universidad de Sevilla, Spain, \email{yhinojos@us.es}\\
\textsuperscript{*} Presenting author. \\
}
{
Facility Location, Network Design
}
\end{abstract_online}


%
%
%

Due to the urgent need to combat climate change, most of the World Organizations have recognized the importance of transitioning to cleaner and more sustainable energy sources in the coming years to reduce greenhouse gas emissions. It is a known fact that one of the main strategies to achieve the proposed objective is the use of biogas as an alternative renewable energy source to carbon-based energies. In this sense, mathematical optimization is recognized as a fundamental tool for the design and modeling of multiple logistics, transportation and supply chain problems and, in particular, for designing robust and efficient supply chains to integrate bioenergy into any economy~\cite{egieya_jafaru_m_biogas_2018,jensen_optimizing_2017}.
The paper by \cite{blanco_2023}  introduces  a new logistic problem for waste-to-biogas transformation. The authors provide a general and flexible mathematical optimization model that allows decision makers to optimally determine the locations of different types of plants (pretreatment, biogas, and liquefaction plants), as well as the most efficient distribution of products (waste to biomethane) along the supply chain involved in the logistic process. They jointly integrate all these complex stages into a mixed integer linear programming (MILP) model that attempts to minimize the overall transportation cost of the system assuming that a limited budget is available to install all types of plants.
In this work, we  further explore the model proposed in \cite{blanco_2023} by analyzing different alternatives to incorporate routes of vehicles in different phases of this problem. Specifically, we provide models to decide the most efficient way to determine the routes to be followed by vehicles passing through all farms (pick-up points) and  pre-treatment plants (delivery points). The goal is to provide a mathematical optimization based framework for choosing the delivery points to open (within a given  set of potential delivery locations) in order to minimize the overall transportation costs required to collect all the demand from the pick-up points and deliver it to one of the open delivery points using a single vehicle. Among other alternatives, we study the problem with the following hypotheses: a single visit to each pick-up point is allowed or multiple visits are allowed; the resulting route is required to be a cycle (it must start and end at the same delivery point) or not; the resulting route must be connected or not (which would imply the use of multiple vehicles, etc.). The different alternatives pose different mathematical challenges for both modeling and solving the resulting problems.
We apply and test the different models in the pickup and delivery problem (PDP)  instances TS2004t2 and TS2004t3 used in \cite{hernandezyJJ2004} and in a real-word dataset based on the region of upper Yahara Watershed in the state of Wisconsin~\cite{SAMPAT2021105199}.

%% file: abstracts/tex/IWOLOCA24_Abstract_1691.tex


\index[author]{Baldassarre, Silvia \\\hspace{0.5cm}
	 University of Naples Federico II, Italy, \email{silvia.baldassarre@unina.it}}
\index[author]{Bruno, Giuseppe \\\hspace{0.5cm}
	 University of Naples Federico II, Italy, \email{giuseppe.bruno@unina.it}}
\index[author]{Cavola, Manuel \\\hspace{0.5cm}
	Pegaso University, Italy,  \email{manuel.cavola@pegaso.it}}
\index[author]{Pipicelli, Eduardo \\\hspace{0.5cm}
	 University of Naples Federico II, Italy, \email{eduardo.pipicelli@unina.it}}


\title{A mathematical model for designing a postal banking network}


\author{
 Silvia Baldassarre\,\textsuperscript{a}, Giuseppe Bruno\,\textsuperscript{a}, Manuel Cavola\,\textsuperscript{b}, Eduardo Pipicelli\,\textsuperscript{a,*}}

\address{
      \textsuperscript{a } University of Naples Federico II, Italy, Department of Industrial Engineering, Piazzale Tecchio 80, Naples, Italy,   \email{silvia.baldassarre@unina.it},   \email{giuseppe.bruno@unina.it},  \email{eduardo.pipicelli@unina.it}\\
      \textsuperscript{b } Pegaso University, \email{manuel.cavola@pegaso.it}\\
      \textsuperscript{*} Presenting author. \\
}

\keywords{Facility Location, Network design, Decision support}

\begin{abstract_online}
{A mathematical model for designing a postal banking network}
{
 Silvia Baldassarre\,\textsuperscript{a}, Giuseppe Bruno\,\textsuperscript{a}, Manuel Cavola\,\textsuperscript{b}, Eduardo Pipicelli\,\textsuperscript{a,*}
}
{
\textsuperscript{a }University of Naples Federico II, Italy, Department of Industrial Engineering, Piazzale Tecchio 80, Naples, Italy,   \email{silvia.baldassarre@unina.it},   \email{giuseppe.bruno@unina.it},  \email{eduardo.pipicelli@unina.it}\\
\textsuperscript{b }Pegaso University, \email{manuel.cavola@pegaso.it}\\
\textsuperscript{*} Presenting author. \\
}
{
Facility Location, Network design, Decision support
}
\end{abstract_online}


Banking groups worldwide are significantly closing their branches while leveraging digital channels to deliver various banking services, from basic to complex ones. 
This shift has raised concerns about financial exclusion, particularly for people hesitant or reluctant to adopt digital channels or those living in lower-income and less densely populated areas.
In this scenario, postal operators have emerged as main competitors in providing physical and face-to-face financial services in areas without a banking presence, leveraging their capillary network of post offices.

In this work, we propose a mathematical model to support postal providers in making decisions about the network design for consulting services. 
The objective is to maximize the demand from customers who are financially excluded due to branch closures. 
The postal banking network design concerns activating consulting services in post offices to capture market share and planning service time availability to satisfy the demand.
An available budget for activating consulting services is considered.  

The model has been applied to a real case study: the rural areas of the Campania region in southern Italy. 
In these areas, bank branch closures have significantly impacted customers' accessibility to banking services over the last few years. 
The results demonstrate that the model can effectively support the consulting network design process, offering an effective framework for making informative decisions based on the preferred postal provider's provision policy. 

%% file: abstracts/tex/IWOLOCA24_Abstract_8661.tex


\index[author]{Baldomero-Naranjo, Marta \\\hspace{0.5cm}
	 Universidad de C\'adiz, Spain, \email{marta.baldomero@uca.es}}
\index[author]{G\'azquez, Ricardo \\\hspace{0.5cm}
	 Universidad Carlos III de Madrid, Spain, \email{ricardo.gazquez@uc3m.es}}
\index[author]{Rodr\'iguez-Ch\'ia, Antonio M. \\\hspace{0.5cm}
	 Universidad de C\'adiz, Spain, \email{antonio.rodriguezchia@uca.es}}


\title{An upgrading approach for the multi-type Maximal Covering Location Problem}


\author{
 Marta Baldomero-Naranjo\textsuperscript{a}, Ricardo G\'azquez\,\textsuperscript{b,*}, Antonio M. \\Rodr\'iguez-Ch\'ia\,\textsuperscript{a}}

\address{
      \textsuperscript{a }Universidad de C\'adiz, Spain,  \email{marta.baldomero@uca.es},   \email{antonio.rodriguezchia@uca.es}\\
      \textsuperscript{b }Universidad Carlos III de Madrid, Spain, \email{ricardo.gazquez@uc3m.es}\\
}


\begin{abstract_online}
{An upgrading approach for the multi-type Maximal Covering Location Problem}
{
 Marta Baldomero-Naranjo\,\textsuperscript{a}, Ricardo G\'azquez\,\textsuperscript{b,*}, Antonio M. \\Rodr\'iguez-Ch\'ia\,\textsuperscript{a}
}
{
      \textsuperscript{a }Universidad de C\'adiz, Spain,  \email{marta.baldomero@uca.es},   \email{antonio.rodriguezchia@uca.es}\\
      \textsuperscript{b }Universidad Carlos III de Madrid, Spain, \email{ricardo.gazquez@uc3m.es}\\
      \textsuperscript{*} Presenting author. \\
}
{
Covering Problems, Upgrading, Maximal Covering
}
\end{abstract_online}


Covering location problems are among core problems in Location Science. These problems arise when facilities offer services within a specific distance or maximum time. Depending on the context, different problems may occur, but they generally fall into two categories: set covering and maximal covering.

This presentation centers on the Maximal Covering Location Problem (MCLP), initially introduced by Church and ReVelle (\cite{ChurchReVelle:1974}). Specifically, it explores an extension of MCLP that includes multiple types of facilities and their upgrades.

While most existing literature addresses problems within a single setting, recent research examines the use of multiple settings within the same problem in MCLP (\cite{blanco2023multi}). Conversely, upgrading problems are becoming increasingly prominent in recent location literature. These extensions involve coordinating facility location decisions with improvements to related infrastructure. In the MCLP, authors in \cite{upMCLP22} examine scenarios where upgrading the network reduces the distance or time from a facility to a customer.

We named the model presented in this talk as the Multi-type Maximal Covering Location Problem with Upgrading (MTMCLP-U). It extends the MCLP by incorporating the two described approaches above: different types of facilities (in both continuous and network settings) and two types of upgrades (for the facilities and the network).

In this model, given a set of customers—some with access to the network and some without—the objective is to place a fixed number of facilities at network nodes and another fixed number in the continuous space to maximize covered demand. Additionally, the model includes two types of binary upgrades (yes-or-no decisions) within a common budget. In the network setting, upgrades reduce the lengths of the edges, while in the continuous setting, upgrades increase the coverage radius, thereby reaching more customers.

\begin{figure}[H]
	\centering
	\begin{subfigure}[b]{.48\linewidth}
		\centering
		\fbox{\includegraphics[scale=0.28]{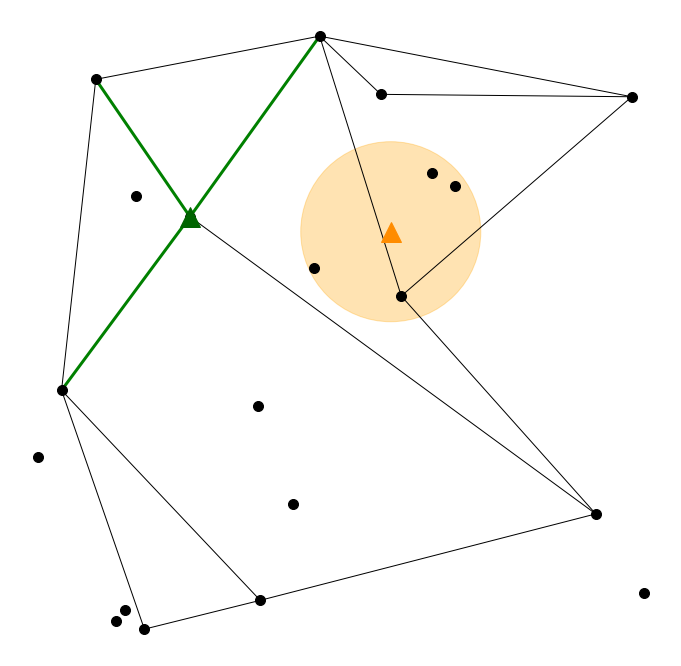}}
		\caption{Without upgrading.} \label{figa}
	\end{subfigure}~
	\begin{subfigure}[b]{.48\linewidth}
		\centering 
		\fbox{\includegraphics[scale=0.28]{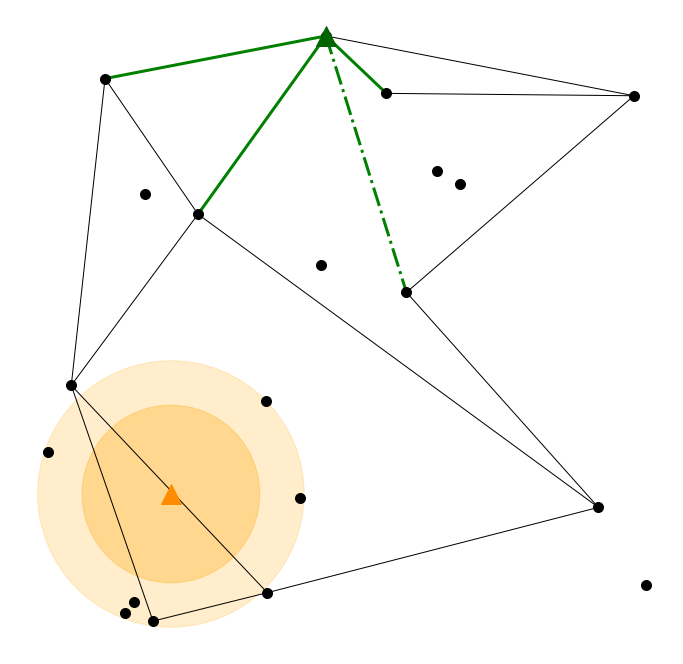}}
		\caption{Upgrading in both settings.} \label{figb}
	\end{subfigure}
	\caption{The green triangle represents a discrete facility, located only at network nodes. The orange triangle represents a continuous facility, which can be located anywhere in the plane. The dashed line indicates an upgraded edge, and light yellow denotes the upgraded coverage area. Figure \ref{figa} shows the facility locations without upgrades. Figure \ref{figb} shows the locations with upgrades, covering more points.}
	\label{fig}
\end{figure}

The MTMCLP-U responds to problems in real applications such as telecommunication networks. For instance, in rural areas where access to telecommunication services cannot be carried over the wired network, repeaters are necessary. See Figure \ref{fig} for a comparison of facility locations with (Fig. \ref{figa}) and without (Fig. \ref{figb}) upgrades.

%% file: abstracts/tex/IWOLOCA24_Abstract_0233.tex

\index[author]{Baldomero-Naranjo, Marta \\\hspace{0.5cm}
	 Universidad de C\'adiz, Spain, \email{marta.baldomero@uca.es}}
\index[author]{Kalcsics, J\"org \\\hspace{0.5cm}
	 The University of Edinburgh, United Kingdom, \email{joerg.kalcsics@ed.ac.uk}}
\index[author]{Rodr\'iguez-Ch\'ia, Antonio M. \\\hspace{0.5cm}
	 Universidad de C\'adiz, Spain, \email{antonio.rodriguezchia@uca.es}}


\begin{abstract_online}
{
The maximal covering location problem with edge downgrades}
{
Marta Baldomero-Naranjo\,\textsuperscript{a,*}, J\"org Kalcsics\,\textsuperscript{b}, Antonio M. Rodr\'iguez-Ch\'ia\,\textsuperscript{a}
}
{
\textsuperscript{a }Universidad de C\'adiz,  \email{marta.baldomero@uca.es},   \email{antonio.rodriguezchia@uca.es}\\
      \textsuperscript{b }The University of Edinburgh, \email{joerg.kalcsics@ed.ac.uk}\\
\textsuperscript{*} Presenting author. \\
}
{
Facility Location, Bilevel optimization, Heuristic algorithms
}
\end{abstract_online}

\section*{Introduction}

This talk deals with an extension of the well-known   Maximal Covering Location Problem (MCLP), see \cite{Church_and_Revele74}. It considers that an agent (attacker) can modify the length of some edges within a budget, increasing the distance from the clients to the facilities. Under this assumption, the Downgrading Maximal Covering Location Problem emerges as a significant challenge, involving the interests of both the facility location planner and potential attackers.  Note that this implies that the distances in the network are modified during the optimization problem, thus, the distance between two nodes after the downgrades will have to be computed within the model. 

This study falls within the domain of downgrading (upgrading) network problems, in which specific elements initially treated as fixed inputs in the classical version are now transformed into decision variables (e.g. the length of the edges). As far as we know, it is the first study analysing edge downgrades. The version of the problem where edge upgrades are considered has been discussed in \cite{BalKalMarRod22, BalKalRod24}. 

\section*{Contributions}
The problem entails the strategic location of facilities within the nodes of a network to maximise coverage while anticipating and mitigating potential attacks aimed at reducing this coverage. Then, it comprises an interaction between two distinct actors, each with conflicting objectives:
\begin{itemize}
    \item The location planner aims to maximise the covered demand while anticipating that an attacker will attempt to reduce coverage by increasing the length of the edges.
    \item The attacker seeks to maximise the demand initially covered by the facilities but left uncovered after the downgrade. To achieve this, they can increase the length of certain edges within a specified budget.
\end{itemize}

We propose a mixed-integer linear bilevel formulation to model the problem. Moreover, we introduce a strategy to preprocess the data to reduce the number of variables and constraints of the formulation. 
Furthermore, a matheuristic algorithm to solve the problem is developed. For this algorithm, we present some variants (levels of intensity), which allow the user to decide the best trade-off between the computational time employed and the quality of the solution.  

Several computational experiments are carried out to illustrate the advantages of applying the introduced model rather than using the classical MCLP and to show the potential and limitations of the proposed matheuristic algorithm. The solutions provided by the matheuristic are compared with the ones obtained by the general bilevel solver developed in \cite{FisLjuMonSin17}.



%% file: abstracts/tex/IWOLOCA24_Abstract_6905.tex


\index[author]{Baldomero-Naranjo, Marta \\\hspace{0.5cm}
	 Universidad de C\'adiz, Spain, \email{marta.baldomero@uca.es}}
\index[author]{Mancuso, Andrea \\\hspace{0.5cm}
	 University of Naples Federico II, Italy, \email{andrea.mancuso@unina.it}}
\index[author]{Masone, Adriano \\\hspace{0.5cm}
	 University of Naples Federico II, Italy, \email{adriano.masone@unina.it}}
\index[author]{Rodr\'iguez-Ch\'ia, Antonio M. \\\hspace{0.5cm}
	 Universidad de C\'adiz, Spain, \email{antonio.rodriguezchia@uca.es}}


\title{New results in the covering tour problem with edge upgrades}


\author{
 Marta Baldomero-Naranjo\,\textsuperscript{a}, Andrea Mancuso\,\textsuperscript{b}, Adriano Masone\,\textsuperscript{b}, Antonio M. Rodr{\'\i}guez-Ch{\'\i}a\,\textsuperscript{a,*}}

\address{
      \textsuperscript{a }Departamento de Estad{\'\i}stica e Investigaci\'on Operativa, Facultad de Ciencias, 
Universidad de C\'adiz, C\'adiz, Spain,  \email{marta.baldomero,antonio.rodriguezchia@uca.es}\\
      \textsuperscript{b }Department of Electrical Engineering and Information Technology, 
University of Naples Federico II, Naples, Italy, \email{andrea.mancuso,adriano.masone@unina.it}\\
}

\keywords{Facility Location, Covering, TSP, upgrading problem}

\begin{abstract_online}
{New results in the covering tour problem with edge upgrades}
{
 Marta Baldomero-Naranjo\,\textsuperscript{a}, Andrea Mancuso\,\textsuperscript{b}, Adriano Masone\,\textsuperscript{b}, Antonio M. Rodr{\'\i}guez-Ch{\'\i}a\,\textsuperscript{a,*}
}
{
\textsuperscript{a }Departamento de Estad{\'\i}stica e Investigaci\'on Operativa, Facultad de Ciencias, Universidad de C\'adiz, C\'adiz, Spain,  \email{marta.baldomero,antonio.rodriguezchia@uca.es}\\
\textsuperscript{b }Department of Electrical Engineering and Information Technology, University of Naples Federico II, Naples, Italy, \email{andrea.mancuso,adriano.masone@unina.it}\\
\textsuperscript{*} Presenting author. \\
}
{
Facility Location, Covering, TSP, Upgrading problem
}
\end{abstract_online}

This talk addresses a version of the Covering Tour Problem,   
a covering problem within the context of location decision
problems  combined with a routing problem. This version considers that  the edges can be upgraded subject to a budget constraint. Therefore, the goal is to find a minimum length tour that cover all the nodes of a graph taking into account that some edges can be upgraded.

\section*{Introduction}

Routing problems are a class of problems focused on finding the most effective way
to move goods, people, or information from one point to another. The primary objective is often to minimize travel time, distance, or operational
costs while adhering to various constraints and requirements.  These problems are central to many industries, including logistics,
supply chain management, and transportation.
Among these routing problems, one of the most fundamental is the Traveling
Salesman Problem (TSP). The TSP involves finding the shortest possible route that
allows a salesman to visit each city in a given list exactly once and return to the
origin city,  see \cite{TSP}.

Coverage problems arise when it is necessary to determine the optimal location
of one or more facilities for the production of goods or the provision of services,
with the goal of meeting the demand for goods or services at various nodes. The
main objective is to minimize the placement costs of facilities or maximize the coverage
of demand points. For each facility, a coverage radius is defined, representing
the distance within which the demand for goods or services can be satisfied.

Combining TSP and convering problems is obtained the covering tour problem (CTP). This problem
 aims to determine a minimum lengh tour over a subset of nodes that
can be visited,  in such a way that the tour contains all these nodes, and every
node  is covered by the tour, see \cite{CTSP}.

The goal of this talk is to analyze the covering tour problem with edge upgrading (CTPEU).
This problem is an extension of the CTP where we consider the possibility of upgrading edges. It combines the
routing aspect seen in the TSP, the covering part of the CTP, and additionally
allows for the possibility of updating edges, similar to the Maximal Covering Problem
with Edge Upgrade, see \cite{Marta}.  Edge upgrading consists in reducing their length (travel time),
typically within certain limits, at a given cost that is proportional to the upgrade.
The total cost of all upgrades is subject to a budget constraint. On this basis, the
CTPEU involves finding a tour of minimum length subject to covering and budget
constraints.


%% file: abstracts/tex/IWOLOCA24_Abstract_5013.tex


\index[author]{Blanco, V\'ictor \\\hspace{0.5cm}
	 Universidad de Granada, Spain, \email{vblanco@ugr.es}}
\index[author]{Gagrani, Praful \\\hspace{0.5cm}
University of Wisconsin-Madison, India,  \email{praful.gagrani@gmail.com}}
\index[author]{Gonz\'alez, Gabriel \\\hspace{0.5cm}
	 Universidad de Granada, Spain, \email{ggdominguez@ugr.es}}


\title{Expanding Reaction Networks through Autocatalytic Subnetworks that Maximize Growth Factor}


\author{
 Blanco, V\'ictor \textsuperscript{a}\, ; Gagrani, Praful \textsuperscript{b}\, ;  Gonz\'alez, Gabriel \textsuperscript{c,*} \, .}

\address{
\textsuperscript{a }Universidad de Granada,  \email{vblanco@ugr.es} \, ,\\
\textsuperscript{b }University of Wisconsin-Madison, \email{praful.gagrani@gmail.com} \, ,\\
\textsuperscript{c }Universidad de Granada, \email{ggdominguez@ugr.es} \, .\\
\textsuperscript{*} Presenting author. \\
}


\begin{abstract_online}
{Expanding Reaction Networks through Autocatalytic Subnetworks that Maximize Growth Factor}
{
V\'ictor Blanco\,\textsuperscript{a}, Praful Gagrani\,\textsuperscript{b},  Gabriel Gonz\'alez\,\textsuperscript{c,*}
}
{
\textsuperscript{a }Universidad de Granada,  \email{vblanco@ugr.es}\\
\textsuperscript{b }University of Wisconsin-Madison \email{praful.gagrani@gmail.com}\\
\textsuperscript{c }Universidad de Granada, \email{ggdominguez@ugr.es}\\
\textsuperscript{*} Presenting author. \\
}
{
Chemical Reaction Networks, Origin Of Life, Autocatalytic Network, Lineal Programming, Network Design
}
\end{abstract_online}


The dynamics of Chemical Reaction Networks (CRNs) and their implications are intricate matters with an important role in the origin of life. This work focuses on identifying autocatalytic subnetworks within CRNs, crucial for explicating self-replication dynamics in biological and prebiotic systems. Autocatalysis, where a set of species catalyze their own production, is pivotal for the emergence of life and chemical evolution. Detecting autocatalytic species within CRNs poses a significant computational challenge, addressed here through a mathematical optimization-based framework. Several tools have been developed to detect and enumerate autocatalytic subnetworks by solving certain mixed integer optimization problems~\cite{peng_hierarchical_2022,gagrani2023geometry}. It is demonstrated that Discrete Optimization is a powerful tool in the search for these complex structures. The proposed model intends to determine the most prominent autocatalytic subnetwork within a CRN and simulate its expansion over time. Through sequential incorporation of species and reactions, the model constructs increasingly autocatalytic subnetworks, shedding light on the dynamics of CRN expansion. The study offers insights into fundamental questions about origins of chemical systems, CRN expansion over time, and the prerequisites for generating autocatalytic species.
\vspace{-0.1cm}

%% file: abstracts/tex/IWOLOCA24_Abstract_4203.tex


\index[author]{Blanco, V\'ictor \\\hspace{0.5cm}
	 Universidad de Granada, Spain, \email{vblanco@ugr.es}}
\index[author]{Jap\'on, Alberto \\\hspace{0.5cm}
	 Universidad de Sevilla, Spain, \email{ajapon1@us.es}}
\index[author]{Puerto, Justo \\\hspace{0.5cm}
	 Universidad de Sevilla, Spain, \email{puerto@us.es}}
\index[author]{Zhang, Peter \\\hspace{0.5cm}
	 Carnegie Mellon University, United States of America, \email{pyzhang@cmu.edu}}


\title{Classification forests via mathematical programming: case study in obesity detection}


\author{
 Víctor Blanco\,\textsuperscript{a} , Alberto Japón\,\textsuperscript{b,*}\ , Justo Puerto\,\textsuperscript{b}, Peter Zhang\,\textsuperscript{c}}

\address{
\textsuperscript{a} Universidad de Granada,  \email{vblanco@ugr.es}\\
\textsuperscript{b} Universidad de Sevilla, \email{ajapon1@us.es}
\email{puerto@us.es}\\
\textsuperscript{c} Carnegie Mellon University,  \email{pyzhang@cmu.edu}\\
\textsuperscript{*} Presenting author. \\
}


\begin{abstract_online}
{Classification forests via mathematical programming: case study in obesity detection}
{
Víctor Blanco\,\textsuperscript{a}, Alberto Japón\,\textsuperscript{b,*}, Justo Puerto\,\textsuperscript{b}, Peter Zhang\,\textsuperscript{c}
}
{
\textsuperscript{a} University of Granada,  \email{vblanco@ugr.es}\\
\textsuperscript{b} University of Seville, \email{ajapon1@us.es}
\email{puerto@us.es}\\
\textsuperscript{c} Carnegie Mellon University,  \email{pyzhang@cmu.edu}\\
\textsuperscript{*} Presenting author. \\
}
{
Classification Forests, Classification Trees, Interpretable Machine Learning
}
\end{abstract_online}


In recent years, there has been an enormous development in the field of machine learning and its applications in both academia and industry. One of the most prominent lines of research has been interpretable machine learning, due to the general preference of practitioners to use predictive models that are understandable from a human point of view over black-box ones \cite{IML}. Whether models are more or less understandable for a given audience may have a certain degree of subjectivity. However, there is agreement that the simpler the model, the better. Mathematical optimization has made great contributions in this respect when working with moderate sample sizes. As an example we can look at the field of classification trees, where optimal classification trees achieve accuracy rates equal to or better than heuristic trees with more splits \cite{OCT}, or at decision forests, where optimal classiication forests  need fewer trees than random forests \cite{OCF}.\\

On the other hand, in many applications, it is not only important to understand the overall classifier solution, but also to understand certain individual cases of sample predictions when the prediction has not been as desired. This is for example the case of a person who is denied a mortgage loan, where it is beneficial for both the individual and the entity to understand what the individual could change to classify as valid for the loan. These are the so-called counterfactual explanations. In this work, we will use classification forests and mathematical programming in a case study applied to the detection of obesity, and we will study the individual solutions of the people in the sample to understand what can be done to prevent obesity.


%% file: abstracts/tex/IWOLOCA24_Abstract_9068.tex


\index[author]{Blanco, V\'ictor \\\hspace{0.5cm}
	 Universidad de Granada, Spain, \email{vblanco@ugr.es}}
\index[author]{Mart\'inez-Ant\'on, Miguel \\\hspace{0.5cm}
	 Universidad de Granada, Spain, \email{mmanton@ugr.com}}


\title{Finding a smallest mediated set}


\author{
 V\'ictor Blanco\,\textsuperscript{a,b}, Miguel Mart\'inez-Ant\'on\,\textsuperscript{a,b,*}}

\address{
      \textsuperscript{a }Institute of Mathematics (IMAG), University of Granada,  \email{vblanco@ugr.com},   \email{mmanton@ugr.com}\\
      \textsuperscript{b }Dpt. Quantitative Methods for Economics \& Business, University of Granada\\
}


\begin{abstract_online}
{Finding a smallest mediated set}
{
  V\'ictor Blanco\,\textsuperscript{a,b}, Miguel Mart\'inez-Ant\'on\,\textsuperscript{a,b,*}
}
{
\textsuperscript{a }Institute of Mathematics (IMAG), University of Granada,  \email{vblanco@ugr.com},   \email{mmanton@ugr.com}\\
      \textsuperscript{b }Dpt. Quantitative Methods for Economics \& Business, University of Granada\\
      \textsuperscript{*} Presenting author. \\
}
{
Discrete Location, Mediated Sets, Sums of Squares, Circuit Polynomials
}
\end{abstract_online}


\section*{Introduction}

Since Hilbert's 17th Problem, which states if every real homogeneous polynomial (form) that is nonnegative can be written as sum of squares (SOS) of suitable forms~\cite{hilbert}, the study of nonnegativity of real, multivariate polynomials and sums of squares is a key problem in real algebraic geometry and an active field of research especially due to their useful applications in polynomial optimization developed in the last two decades. 

On the one hand, a circuit polynomial is a polynomial $p\in \mathbb{R}[x_1,\ldots, x_d]$ of the form
$$
p(\textbf{x})=\displaystyle\sum_{j=0}^rc_{\boldsymbol{\alpha}(j)}\textbf{x}^{\boldsymbol{\alpha}(j)}+c_{\boldsymbol{\beta}}\textbf{x}^{\boldsymbol{\beta}}
$$
where $r\leq d$, exponents $\boldsymbol{\alpha}(j),\boldsymbol{\beta}\in \mathbb{Z}^d_+$, and coefficients $c_{\boldsymbol{\alpha}(j)}>0$ and $c_{\boldsymbol{\beta}}\in \mathbb{R}^*$ such that $\mathrm{New}(p)$ is a simplex with even vertices $\boldsymbol{\alpha}(0),\ldots, \boldsymbol{\alpha}(r)$ and the exponent $\boldsymbol{\beta}$ is in the strict interior of $\mathrm{New}(p)$.

On the other hand, a mediated set $L$ of an integral simplex $S$ with vertex, $\mathrm{Vert}(S)$, in $(2\mathbb{Z})^d$ is a subset of lattice points in $\mathbb{Z}^d\cap S$ satisfying:
\begin{itemize}
    \item $\mathrm{Vert}(S$)$\subset L$, and

    \item if $\boldsymbol{\gamma}_i\in L$, then there exist $\boldsymbol{\gamma}_j$,$\boldsymbol{\gamma}_k\in (2\mathbb{Z})^d\cap L$ with $\boldsymbol{\gamma}_i=\frac{1}{2}(\boldsymbol{\gamma}_j+\boldsymbol{\gamma}_k).$
\end{itemize}

Iliman and de Wolff~\cite{iliman} proved that a nonnegative circuit polynomial is SOS if and only if there exists a mediated set $M$ of its Newton polytope $S$ that contains the exponent of one distinguished term ($\boldsymbol{\beta}$) and, in this case, it will be a linear combination with strictly positive coordinates of
$$
\left\{\left(\textbf{x}^{\frac{\boldsymbol{\gamma}_j}{2}}-\textbf{x}^{\frac{\boldsymbol{\gamma}_k}{2}}\right)^2\right\}_{\boldsymbol{\gamma}_i\in M\setminus \mathrm{Vert}(S).}
$$

\begin{figure}[H]
\begin{center}
\fbox{\includegraphics[scale=0.5]{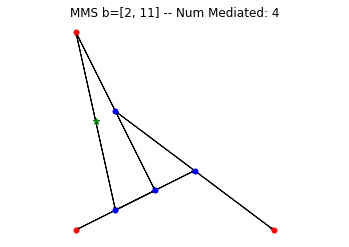}}
\end{center}
\caption*{Minimum mediated set (blue dots) for $\mathrm{Vert}(S)=\{(0,0),(20,0),(0,20)\}$ (red dots), and the interior point $(2,11)$ (green star). Solution associated to $2x^{20}+11y^{20}+7-20x^2y^{11}$.}
\end{figure}
In~\cite{blanco}, the authors study mediated sets in the continuous case due to their crucial role in the minimal second-order cone representation of rational generalized power cones. They obtain results on the combinatorial structure of these sets using mathematical optimization tools. This approach is applied both in the continuous case to study generalized power cones and in the discrete case to address a question arising from Hilbert's 17th Problem. Given a simplex $S$ with even vertices and a strictly interior point $\boldsymbol{\beta}$, what is the minimum mediated set of $S$ that contains $\boldsymbol{\beta}$? In case this set is empty we have a certification of a nonnegative circuit polynomial supported on $\mathrm{Vert}(S)\cup \{\boldsymbol{\beta}\}$ is not SOS, otherwise, we will know a minimal representation of that polynomial as SOS, and the involved binomial.

The authors' first approach requires enumerating every single even point in the simplex $S$, which presents a significant challenge. The next step involves finding a small set of potential mediated points to determine the minimum mediated set. The authors propose a constructive method to build this potential set, starting from the vertex set of the simplex and the interior point. They also introduce suitable mathematical optimization models to find the optimal location of the mediated points in the mediated set while enumerating as few points as possible.


%% file: abstracts/tex/IWOLOCA24_Abstract_7341.tex


\index[author]{Blanco, V\'ictor \\\hspace{0.5cm}
	 Universidad de Granada, Spain, \email{vblanco@ugr.es}}
\index[author]{Pozo, Miguel A. \\\hspace{0.5cm}
	 Universidad de Sevilla, Spain, \email{pozo@us.es}}
\index[author]{Puerto, Justo \\\hspace{0.5cm}
	 Universidad de Sevilla, Spain, \email{puerto@us.es}}
\index[author]{Torrej\'on, Alberto \\\hspace{0.5cm}
	 Universidad de Sevilla, Spain, \email{atorrejon@us.es}}


\title{Modeling ordered interactions in Location Problems}


\author{
 Víctor Blanco\,\textsuperscript{a}, 
 Miguel A. Pozo\,\textsuperscript{b}, 
 Justo Puerto\,\textsuperscript{b}, 
 Alberto Torrejón\,\textsuperscript{b,*}}

\address{
      \textsuperscript{a }Department of Quantitative Methods for Economics \& Business, University of Granada, Spain. \\ Institute of Mathematics of the University of Granada (IMAG), Spain.  \\ \email{vblanco@ugr.es},   \\
      \textsuperscript{b }Department of Statistics and Operational Research, University of Seville, Spain. \\ Institute of Mathematics of the University of Seville (IMUS), Spain. \\ \email{pozo@us.es}, \email{puerto@us.es}, \email{atorrejon@us.es}\\
      \textsuperscript{*} Presenting author. \\
}


\begin{abstract_online}
{Modeling ordered interactions in Location Problems}
{
 Víctor Blanco\,\textsuperscript{a}, 
 Miguel A. Pozo\,\textsuperscript{b}, 
 Justo Puerto\,\textsuperscript{b}, 
 Alberto Torrejón\,\textsuperscript{b,*}
}
{
\textsuperscript{a }Department of Quantitative Methods for Economics \& Business, University of Granada, Spain. \\ Institute of Mathematics of the University of Granada (IMAG), Spain,  \\ \email{vblanco@ugr.es},   \\
\textsuperscript{b }Department of Statistics and Operational Research, University of Seville, Spain. \\ Institute of Mathematics of the University of Seville (IMUS), Spain, \\ \email{pozo@us.es}, \email{puerto@us.es}, \email{atorrejon@us.es}\\
\textsuperscript{*} Presenting author. \\
}
{
Discrete location, Ordered optimization, Quadratic optimization
}
\end{abstract_online}

\section*{Introduction}

Ordered optimization (see \cite{ref1} \& \cite{ref2}) allows a straightforward and flexible generalization of facility location problems, since many well-known problems in the literature can be modelled using an ordered formulation. 
The ordered description of these problems relies on the vector of weights, usually named $\lambda$-vector, that multiplies the ordered costs in the objective function.
By means of an appropriate $\lambda$-vector, many location problems such as median, center or centdian problems can be modelled, but also, since negative and non-monotonous weight vectors are allowed, one can also consider obnoxious location problems, preference modelling or other abstract concepts such as, for example, fairness. 

However, the range of problems that can be studied within an ordered framework is limited to the linearity of the objective function.
By using a quadratic approach to ordered optimization, and considering a matrix of weights instead of a vector, the number of criteria to be modelled is considerably multiplied, allowing, among others, to minimize interaction between costs or to model more specific objective functions which require quadratic terms, such as the variance. 
On the other hand, introducing ordered interactions as a modeling feature increases the complexity of the problems at hand, which motivates the study of efficient resolution methods that allow the scalability of these problems.

This work motivates this new approach, which allows a higher level of generalisation of several combinatorial problems, particularly location problems, providing a mathematical formalization of this approach, describing exact models as solution methods and performing empirical comparison between these different solutions methods.


%% file: abstracts/tex/IWOLOCA24_Abstract_5262.tex

\index[author]{B\"using, Christina \\\hspace{0.5cm}
	 RWTH Aachen University, Germany, \email{buesing@combi.rwth-aachen.de}}
\index[author]{Leitner, Markus \\\hspace{0.5cm}
	 Vrije Universiteit Amsterdam, The Netherlands, \email{m.leitner@vu.nl}}
\index[author]{Wrede, Sophia \\\hspace{0.5cm}
	 RWTH Aachen University, Germany, \email{wrede@combi.rwth-aachen.de}}

\begin{abstract_online}
{Cover-based inequalities for the single-source capacitated facility location problem with customer preferences
}
{
Christina B\"using\,\textsuperscript{a}, Markus Leitner\,\textsuperscript{b}, Sophia Wrede\,\textsuperscript{a,*}
}
{
\textsuperscript{a }RWTH Aachen University, Aachen, Germany, \email{buesing@combi.rwth-aachen.de}, \\ \email{wrede@combi.rwth-aachen.de} \\
\textsuperscript{b }Vrije Universiteit Amsterdam, Amsterdam, The Netherlands, \email{m.leitner@vu.nl}\\
\textsuperscript{*} Presenting author. \\
}
{
Facility location, preference constraints, valid inequalities
}
\end{abstract_online}


In the classical \textit{single-source capacitated facility location problem} (CFLP), a set of facilities needs to be chosen in order to cover the demand of customers. 
Customers are assigned to any open facility such that the capacity of the facility is not exceeded and the total cost consisting of opening and assignment costs is minimised. 
However, in many real-world applications customers are not willing to travel to any open facility assigned to them but want to select an open facility according to their preferences\,\cite{KANG2023,Polino2023AFL}. 
Such deviations can turn feasible solutions for the CFLP infeasible.
The \textit{single-source capacitated facility location problem with customer preferences} (CFLP-CP) takes this behavior into account by assigning each customer to their most preferred open facility. 

Both the CFLP and CFLP-CP are strongly NP-hard. Preference constraints, however, imply certain combinatorial structures which do not occur in the classical CFLP. 
These implied structures render some previously NP-hard special cases of the CFLP easy to solve\,\cite{BUESING2022}. 
One remarkable structure implied by preference constraints intertwines the assignments of customers with related preferences to the same facility. 

In this talk, we focus on cover-based inequalities for the CFLP-CP. 
We first strengthen the well-known cover inequalities so that they incorporate the specific structure of the CFLP-CP mentioned above.
Afterwards, we strengthen these inequalities by utilising further properties of solutions for the CFLP-CP arising from the combination of capacities and preference constraints.
Lastly, we evaluate the impact of these considered inequalities for two preference types in a computational study.

%% file: abstracts/tex/IWOLOCA24_Abstract_2207.tex


\index[author]{Canca, David \\\hspace{0.5cm}
	 Universidad de Sevilla, Spain, \email{dco@us.es}}
\index[author]{Wang, Xingrong \\\hspace{0.5cm}
	 School of Systems Science, Beijing , China \email{20114009@bjtu.edu.cn}}


\title{A bi-level metaheuristic for the design of hierarchical transit networks in a multimodal context}


\author{
 David Canca\,\textsuperscript{a,*}, Xingrong Wang\,\textsuperscript{b}}

\address{
      \textsuperscript{a }Universidad de Sevilla,  \email{dco@us.es}\\
      \textsuperscript{b }School of Systems Science, Beijing, \email{20114009@bjtu.edu.cn}\\
}


\begin{abstract_online}
{A bi-level metaheuristic for the design of hierarchical transit networks in a multimodal context}
{
 David Canca\,\textsuperscript{a,*}, Xingrong Wang\,\textsuperscript{b}
}
{
\textsuperscript{a }Universidad de Sevilla,  \email{dco@us.es}\\
\textsuperscript{b }School of Systems Science, Beijing, \email{20114009@bjtu.edu.cn}\\
\textsuperscript{*} Presenting author. \\
}
{
Network Design, Hierarchical transportation network, ALNS
}
\end{abstract_online}


\section*{Introduction}

In many large cities around the world, rapid rail rapid transit systems are used as the main transportation mode to meet the mobility needs of citizens. Despite the large volume of passengers transported, rapid transit lines rarely can meet the total demand for passengers. The construction of a well-coordinated metro-bus bimodal system may benefit not only passengers performing bimodal trips in the city area, those who see their travel possibilities increased, but also passengers using each individual transit mode, since sharing the total demand, the congestion of both modes can be easily managed. If both the metro and bus regard themselves as a self-standing system, an uncoordinated metro-bus bimodal configuration with excessive intermodal competition or defective intermodal connection may occur, just like the dilemma faced by the metro and bus systems of some megacities in China. To ensure a safe operation, it has been necessary to implement different passenger demand management measures, such as passenger inflow control systems (\cite{Wang2020, Zhang2021, Hu2023}. However, the bus demand has decreased year after year. Many bus operators have to rely on unsustainable government financial subsidies to maintain daily operations \cite{Liang2019}. Since modifying the bus network is relatively easy and cheaper than acting over the railway network, renewing the bus network design to adapt it to the existing metro system is the most efficient way to improve the bi-modal network functionality. Although several previous studies investigated the effect of the metro on a bus network design taking into account the modal choice for passengers (e.g., \cite{Liang2019, Wan2009}),  no specific attention has been paid to the potential effect of considering different types of bus lines, forming a hierarchical system, working together with an existing metro network. To compensate for the deficiencies existing in both actual operation and theoretical research, this research deals with the problem of simultaneously determining the configuration of a hierarchical bus network and the frequency of the lines to ensure a coordinated operation with an existing metro, explicitly considering the multimodal route choice of passengers.  The hierarchical bus network structure consists of different line types: regular, trunk, and feeder lines. Each type corresponds to a specific stopping pattern -characterized by a particular inter-stop distance, commercial speed, and terminal stations-. The problem is formulated as a non-convex non-linear optimization model. To efficiently solve the model, a bi-level heuristic that breaks down the integrated problem into two interlaced subproblems and then solves them iteratively is presented. In the upper level, an adaptive large neighbourhood search metaheuristic is responsible for proposing new network designs, whereas, in the lower level, a heuristic procedure attempts to determine the best frequencies and perform the passenger assignment. The obtained results demonstrate the capabilities of the proposed approach. As main conclusion, by achieving a trade-off between lines with different running speeds resulting from different average inter-stop distances, hierarchical bus network designs can incentivize commuters into trip plans that contribute to improving the overall integrated metro-bus bimodal network performance.

%% file: abstracts/tex/IWOLOCA24_Abstract_9914.tex


\index[author]{Cattaruzza, Diego \\\hspace{0.5cm}
  Centrale Lille, France, \email{name1@institution1.org}}
\index[author]{Labb\'e, Martine \\\hspace{0.5cm}
  Universit\'e Libre de Bruxelles, Belgium, \email{martine.labbe@ulb.be}}
\index[author]{Petris, Matteo \\\hspace{0.5cm}
  INRIA, France, \email{matteo.petris@inria.fr}}
\index[author]{Roland, Marius \\\hspace{0.5cm}
  Polytechnique Montreal, Canada, \email{mmmroland@gmail.com}}
\index[author]{Schmidt, Martin \\\hspace{0.5cm}
  Trier University, Germany, \email{martin.schmidt@uni-trier.de}}

\title{Novel Valid Inequalities for Chance-\\Constrained Problems
  with Finite Support}


\author{Diego Cattaruzza\,\textsuperscript{a,b},
  Martine Labb\'e\,\textsuperscript{c,b,*},
  Matteo Petris\,\textsuperscript{b},
  Marius Roland\,\textsuperscript{d},
  Martin Schmidt\,\textsuperscript{e}}

\address{
  \textsuperscript{a }Centrale Lille, \email{diego.cattaruzza@centralelille.fr}\\
  \textsuperscript{b }INOCS Team, INRIA, Lille, \email{matteo.petris@inria.fr}\\
  \textsuperscript{c }Universit\'e Libre de Bruxelles, \email{martine.labbe@ulb.be}\\
  \textsuperscript{d }Polytechnique Montreal, \email{mmmroland@gmail.com}\\
  \textsuperscript{e }Trier University, \email{martin.schmidt@uni-trier.de}\\
}


\begin{abstract_online}
{Novel Valid Inequalities for Chance-Constrained Problems
  with Finite Support}
{
Diego Cattaruzza\,\textsuperscript{a,b},
  Martine Labb\'e\,\textsuperscript{c,b,*},
  Matteo Petris\,\textsuperscript{b},
  Marius Roland\,\textsuperscript{d},
  Martin Schmidt\,\textsuperscript{e}
}
{
  \textsuperscript{a }Centrale Lille, \email{diego.cattaruzza@centralelille.fr}\\
  \textsuperscript{b }INOCS Team, INRIA, Lille, \email{matteo.petris@inria.fr}\\
  \textsuperscript{c }Universit\'e Libre de Bruxelles, \email{martine.labbe@ulb.be}\\
  \textsuperscript{d }Polytechnique Montreal, \email{mmmroland@gmail.com}\\
  \textsuperscript{e }Trier University, \email{martin.schmidt@uni-trier.de}\\
  \textsuperscript{*} Presenting author. \\
}
{
Chance constraints, MILP, Valid inequalities
}
\end{abstract_online}


\section*{Introduction}

We consider the (mixed-integer) linear chance-constrained problem
\begin{subequations}
  \label{eq:chanceconstrained}
  \begin{align}
    \min_x \quad
    & c^\top x \\
    \text{s.t.} \quad
    & x \in X, \\
    &	\mathbb{P} ( A^{\omega} x - b^{\omega} \geq 0) \geq 1- \tau
  \end{align}
\end{subequations}

where $X \subset \mathbb{R}^n$ is a non-empty and compact set defined
by deterministic constraints, which possibly include integrality
restrictions on some or all of the variables.
The cost vector is given by $c \in \mathbb{R}^n$.
Let $(\Omega, 2^{\Omega}, \mathbb{P})$ be a discrete and finite
probability space such that $\Omega = \{\omega^s:s \in \mathcal{S}\}$
and $p^s = \mathbb{P}(\omega = \omega^s)$ for $s \in \mathcal{S}$ and
$\sum_{s\in \mathcal{S}}p^s = 1$ holds.
Moreover, let $A^{\omega} \in\mathbb{R}^{m \times n}$ be a matrix of random
constraint coefficients and let $b^{\omega}\in\mathbb{R}^m$ be a
vector of random right-hand sides.
Finally, $\tau \in (0,1)$ is the risk tolerance.

It is well known that this problem can be reformulated as a
Mixed-Integer Linear Problem (MILP) by introducing one binary variable
$y^s$ for each scenario such that $y^s$ takes value $1$ if the
constraint set of that scenario is satisfied by the solution.

We derive two novel families of valid inequalities.
One is based on strong duality of the linear programming formulation
of the respective quantile, whereas the other exploits a covering
argument.
We prove that the covering-based inequalities are stronger than the
first family as well as that they dominate the covering inequalities
proposed in~\cite{Qiu_et_al:2024} for the case of a single constraint
per scenario.
Moreover, we state and prove results on the closure of these valid
inequalities and compare the closure of our second family of
inequalities with the closure of the so-called mixing-set inequalities
introduced in~\cite{Luedtke:2014}.

Finally, we illustrate the impact of our novel valid inequalities in a
numerical study based on instances proposed in~\cite{Song_et_al:2014}.

%% file: abstracts/tex/IWOLOCA24_Abstract_5122.tex


\index[author]{Dom\'inguez, Concepci\'on \\\hspace{0.5cm}
	 University of Murcia, Spain, \email{concepcion.dominguez@um.es}}


\title{Stable Solutions for the Capacitated \\ Simple Plant Location Problem with Order}


\author{
 Concepción Domínguez\,\textsuperscript{a,*} 
 }

\address{
      \textsuperscript{a }University of Murcia,  \email{concepcion.dominguez@um.es},  
     \vspace{.5cm}
}

\keywords{Location with preferences, Simple Plant Location Problem, Facility Location, Combinatorial Optimization, Matchings under preferences, IWOLOCA}

\begin{abstract_online}
{Stable Solutions for the Capacitated Simple Plant Location Problem with Order}
{
Concepción Domínguez\,\textsuperscript{a,*}
}
{
\textsuperscript{a }University of Murcia,  \email{concepcion.dominguez@um.es}\\
\textsuperscript{*} Presenting author. \\
}
{
Location with preferences, Simple Plant Location Problem, Facility Location, Combinatorial Optimization, Matchings under preferences
}
\end{abstract_online}


\section*{Introduction}

The Simple Plant Location Problem (SPLP) is a well-known problem in the Location Field where the aim is to open a set of facilities and allocate each customer to a facility minimizing the total cost of location plus allocation. In the SPLP, it is assumed that the customers have no decision in the allocation, so they are allocated to the open plant that minimizes the cost. The SPLP with Order (SPLPO) arises when the preferences of the customers regarding their allocation is taken into account in the allocation process. Thus, in the SPLPO  customers rank the facilities according to their preferences and they are allocated to their most-preferred facility among the open ones. The SPLPO was introduced in 1987 \cite{Hanjoul1987} and has since been thoroughly studied in the literature (see.\ e.g.\ \cite{Cabezas2023, Canovas2007}). 

In the Capacitated version of the problem (CSPLPO) it is assumed that there is a limited number of customers that can be allocated to each plant. When capacities in the plants and customers' preferences are involved at a time, a decision rule has to be taken to discern which customers are allocated to plants with high demand. As might be expected, any solution with full plants can have \emph{envy customers}, i.e.\ customers that would have liked to be allocated to a different plant which is full. Different decision rules on how to deal with the allocation result in different versions of the problem. Very recently, a bilevel problem has been proposed in \cite{Calvete2020} where the ranking of each customer is given as a list of numerical preferences (the smaller the value, the greater the preference towards the plant) and the aim of the follower problem is to maximize the preferences of the customers \emph{globally}. 

In our work, we propose a new approach where only a \emph{strict total order} of the plants (with no numerical values) is taken into consideration for each customer. In this setting, new decision criteria for the allocation are defined based on the concept of stability in an allocation. New formulations are proposed for each decision criteria, along with valid inequalities and theoretical results designed to obtain compact and tighten them. Preliminary computational results showing the efficiency of the approach and the resultant stable solutions are provided.  

%% file: abstracts/tex/IWOLOCA24_Abstract_2696.tex


\index[author]{Escudero, Laureano F.\\\hspace{0.5cm}
          Universidad Rey Juan Carlos, Spain,  \email{laureano.escudero@urjc.es}}
\index[author]{Gar\'in, Mar\'ia Araceli\\\hspace{0.5cm}
        Universidad del Pa\'is Vasco, Spain, \email{mariaaraceli.garin@ehu.eus}}
\index[author]{Unzueta, Aitziber\\\hspace{0.5cm}
         Universidad del Pa\'is Vasco, Spain, \email{aitziber.unzueta@ehu.eus}}


\title{Cross-docking platforms design and \newline Distributionally robust optimization}


\author{
  Laureano F. Escudero\,\textsuperscript{a,*},
  María Araceli Gar{\'i}n\,\textsuperscript{b}, Aitziber Unzueta\,\textsuperscript{c}}

\address{
      \textsuperscript{a*} Statistics and Operations Research Area, Universidad Rey Juan Carlos
(URJC), Móstoles (Madrid), Spain\email{: laureano.escudero@urjc.es}\\
      \textsuperscript{b} Quantitative Methods Dep., Universidad del Pa\'is Vasco
(UPV/EHU), Bilbao (Bizkaia), Spain\email{:  mariaaraceli.garin@ehu.eus}\\
      \textsuperscript{c} Applied Mathematics Dep., Universidad del Pa\'is Vasco
(UPV/EHU), Bilbao (Bizkaia), Spain\email{: aitziber.unzueta@ehu.eus}\\
}


\begin{abstract_online}
{Cross-docking platforms design and distributionally robust optimization}
{
Laureano F. Escudero\,\textsuperscript{a,*},
María Araceli Gar{\'i}n\,\textsuperscript{b}, Aitziber Unzueta\,\textsuperscript{c}
}
{
\textsuperscript{a} Statistics and Operations Research Area, Universidad Rey Juan Carlos (URJC), Móstoles (Madrid), Spain, \email{laureano.escudero@urjc.es}\\
\textsuperscript{b} Quantitative Methods Dep., Universidad del Pa\'is Vasco (UPV/EHU), Bilbao (Bizkaia), Spain, \email{mariaaraceli.garin@ehu.eus}\\
\textsuperscript{c} Applied Mathematics Dep., Universidad del Pa\'is Vasco (UPV/EHU), Bilbao (Bizkaia), Spain, \email{aitziber.unzueta@ehu.eus}\\
\textsuperscript{*} Presenting author. \\
}
{
Cross-dock door design, two-stage stochastic quadratic
combinatorial optimization, Distributionally robust optimization, Ambiguity sets
}
\end{abstract_online}


The Cross-dock Door Design Problem (CDDP) consists of deciding on the number and capacity of inbound and outbound doors,
minimizing the construction and exploitation cost of the infrastructure. The uncertainty, realized in scenarios, lies in the occurrence of these nodes, the delivering material volume and cost, and the capacity’s disruption of the doors.
Introducing a scheme for generating the ambiguity sets for the second stage uncertain parameters in CDDP.

\section*{Introduction}
Distributionally robust optimization (DRO) is motivated as a
counterpart of the usually unknown underlying probability distribution
(PD)followed by the uncertainty in dynamic problems.
This wok presents a DRO scheme for generating an ambiguity set to
represent the uncertainty in a two-stage scenario tree for the highly
combinatorial CDDP to decide the number and nominal capacity of the
strip and stack doors.

It is assumed the availability of a Nominal Distribution (ND) for the
realization of the uncertain parameters in the second stage
of the tree. That ambiguity set is  generated by considering the
projections of appropriate perturbations of the cumulative
distribution functions of the ND realizations in a set of
modeler-driven PDs (say Normal, Weybull, Gamma and Lognormal
distributions).
 Those perturbations are selected based on the minimization of the
 Wasserstein distance between the parameters’ realizations in  each
 member in the ambiguity set and the ND realization,
 so that a given radius is satisfied.

 The underline assumption is that the realizations of the uncertain
 parameters in modeler-driven groups are identically distributed
 random variables, being independent from the parameters that belong to other groups.

%% file: abstracts/tex/IWOLOCA24_Abstract_2643.tex


\index[author]{Espejo, Inmaculada \\\hspace{0.5cm}
	 Universidad de C\'adiz, Spain, \email{inmaculada.espejo@uca.es}}
\index[author]{Mar\'in, Alfredo \\\hspace{0.5cm}
	 Universidad de Murcia, Spain, \email{amarin@um.es}}
\index[author]{Muñoz-Ocaña, Juan Manuel \\\hspace{0.5cm}
	 Universidad de C\'adiz, Spain, \email{juanmanuel.munoz@uca.es}}
\index[author]{P\'aez, Ra\'ul \\\hspace{0.5cm}
	 Universidad de C\'adiz, Spain, \email{raul.paez@uca.es}}
\index[author]{Rodr\'iguez-Ch\'ia, Antonio M. \\\hspace{0.5cm}
	 Universidad de C\'adiz, Spain, \email{antonio.rodriguezchia@uca.es}}


\title{A further research on Stochastic Single-Allocation Hub Location Problem}


\author{
 Inmaculada. Espejo\,\textsuperscript{a,*}, Alfredo Mar\'in\,\textsuperscript{b}, Juan M. Muñoz-Ocaña\,\textsuperscript{a}, Ra\'ul P\'aez\,\textsuperscript{a}, Antonio M. Rodr\'iguez-Ch\'ia\,\textsuperscript{a}}

\address{
      \textsuperscript{a }Departamento de Estad\'istica e Investigaci\'on Operativa, Universidad de C\'adiz, Spain,  \email{inmaculada.espejo@uca.es}, \email{juanmanuel.munoz@uca.es}, \email{raul.paez@uca.es},\\ \email{antonio.rodriguezchia@uca.es} \\
     \textsuperscript{b }Departamento de Estad\'istica e Investigaci\'on Operativa, Universidad de Murcia, Spain,  \email{amarin@um.es}\\
}


\begin{abstract_online}
{A further research on Stochastic Single-Allocation Hub Location Problem}
{
 Inmaculada Espejo\,\textsuperscript{a,*}, Alfredo Mar\'in\,\textsuperscript{b}, Juan M. Muñoz-Ocaña\,\textsuperscript{a}, Ra\'ul P\'aez\,\textsuperscript{a}, Antonio M. Rodr\'iguez-Ch\'ia\,\textsuperscript{a}
}
{
\textsuperscript{a }Departamento de Estad\'istica e Investigaci\'on Operativa, Universidad de C\'adiz, Spain,  \email{inmaculada.espejo@uca.es}, \email{juanmanuel.munoz@uca.es}, \email{raul.paez@uca.es},\\ \email{antonio.rodriguezchia@uca.es} \\
\textsuperscript{b }Departamento de Estad\'istica e Investigaci\'on Operativa, Universidad de Murcia, Spain,  \email{amarin@um.es}\\
\textsuperscript{*} Presenting author. \\
}
{
Stochastic, Hub, Branch-and-cut
}
\end{abstract_online}


\section*{Introduction}

The Hub Location Problem (HLP) has been addressed by a wide community of operations researchers due to its practical relevance. Many reviews about location problems show the extensive activity in this field and the applications of these problems, see \cite{Alumur2021}, \cite{Campbell2012} and \cite{Contreras2019}, among others.

Most studies on the HLPs are focused on models where the input parameters are assumed to be
fixed and known. However, in most real-life problems, parameters such as demands,
transportation costs, capacity of hubs, et cetera are subject to
variation due to a variety of factors such as population shifts, economic, environmental and political circumstances, quality of the provided services, et cetera.
Hence, some information required for planning is not available a priori, and the associated uncertainties are only resolved once the system is constructed and the
hubs are installed.

This work deals with a Single-Allocation Hub Location Problem (\textrm{SAHLP}) where the amount of product sent between
origins and destinations, the transportation costs and the capacities of the hubs, are uncertain and modelled using random variables with realization only after the hubs are selected.
Different scenarios are considered (with their probabilities) and two decisions have to be made: i) the location of the hubs (this location does not change
from one scenario to another), and ii) the allocation of each origin/destination to these hubs in order to minimize the expected overall costs of the system.
This is a realistic practice because the hubs are located before knowing the real scenario and the allocations are determined when the actual parameters are realized.
This is motivated by the necessity of the network operators  to quickly react to the changes in overall system performance  by adjusting the assignments of the origins/destinations to the hubs when the uncertainty is realized.

A compact integer programming formulation is proposed for the stochastic \textrm{SAHLP} with three variants: uncapacitated, capacitated and $p$-hub problem.
Valid inequalities are developed to reinforce the formulation and added in a branch-and-cut procedure.


%% file: abstracts/tex/IWOLOCA24_Abstract_2489.tex


\index[author]{Camacho-Vallejo, Jos\'e-Fernando \\\hspace{0.5cm}
	 Tecnologico de Monterrey, Mexico, \email{fernando.camacho@tec.mx}}
\index[author]{Fern\'andez-Guti\'errez, Juan Pablo \\\hspace{0.5cm}
	 Universidad de Medell\'in, Colombia, \email{jpfernandez@udemedellin.edu.co}}
\index[author]{Villegas, Juan G. \\\hspace{0.5cm}
	 Universidad de Antioquia, Colombia, \email{juan.villegas@udea.edu.co}}


\title{A Reformulation for the Medianoid Problem with Multipurpose Trips}


\author{
Juan Pablo Fern\'andez-Guti\'errez,\textsuperscript{a}, Juan G. Villegas,\textsuperscript{b},
 Jos\'e-Fernando Camacho-Vallejo,\textsuperscript{c,*}}

\address{
      \textsuperscript{a }Faculty of Basic Science, Universidad de Medell\'in, Colombia, \email{jpfernandez@udemedellin.edu.co}\\
      \textsuperscript{b } ALIADO-Analytics and Research for Decision Making, Department of Industrial Engineering, Universidad de Antioquia, Colombia,
      \email{juan.villegas@udea.edu.co}\\
      \textsuperscript{c }Tecnologico de Monterrey, Escuela de Ingenier\'ia y Ciencias, Mexico,  \email{fernando.camacho@tec.mx}\\
      \textsuperscript{*} Presenting Author.
}


\keywords{Competitive Facility Location, Bilevel Optimization, Medianoid Problem, Travelling Purchaser Problem}

\begin{abstract_online}
{A Reformulation for the Medianoid Problem with Multipurpose Trips}
{
Juan Pablo Fern\'andez-Guti\'errez\,\textsuperscript{a}, Juan G. Villegas\,\textsuperscript{b},
 Jos\'e-Fernando Camacho-Vallejo\,\textsuperscript{c,*}
}
{
\textsuperscript{a }Faculty of Basic Science, Universidad de Medell\'in, Colombia,\\ \email{jpfernandez@udemedellin.edu.co}\\
      \textsuperscript{b }ALIADO-Analytics and Research for Decision Making, Department of Industrial Engineering, Universidad de Antioquia, Colombia,
      \email{juan.villegas@udea.edu.co}\\
      \textsuperscript{c }Tecnologico de Monterrey, Escuela de Ingenier\'ia y Ciencias, Mexico,\\  \email{fernando.camacho@tec.mx}\\
      \textsuperscript{*} Presenting author. \\
}
{
Competitive Facility Location, Bilevel Optimization, Medianoid Problem, Travelling Purchaser Problem
}
\end{abstract_online}


\section*{Introduction}

Competitive Facility Location Problems (CFLPs) involve facilities from different companies offering similar commodities or services. Typically, it is assumed that some facilities are already operational in the market, and new ones will be established to increase or maintain companies' market share \cite{mishra2022location}. Due to its similarity to Stackelberg games, CFLPs often employ leader-follower terminology. Within this hierarchical framework, the leader initiates the process by locating their new facilities. Subsequently, taking the leader's decisions into account, the follower determines the locations of their own facilities.\\

These decisions influence the allocation of customers, which often follows a binary assignment rule. Typically, either the closest distance or the cheapest assignment criteria are used. In our study, we assume that the leader already operates some facilities in the market. Therefore, we are examining the location problem from the follower's perspective, with the aim of maximizing its market share. \\

\section*{Motivation}
The motivation for this study stems from the observation that much of the existing research on CFLPs assumes single-purpose trips. In such scenarios, customers typically make back-and-forth trips to facilities to pick up their demanded commodities. However, when multipurpose trips are considered, customers strategically plan their trips to purchase commodities from different locations.\\

Incorporating routing decisions into the facility location decision-making process is crucial to prevent suboptimal solutions in practice \cite{salhi1989effect}.\\

\section*{The Medianoid Problem with Multipurpose Trips}

In this research, we address the Medianoid Problem with Multipurpose Trips (MPMT) as proposed in \cite{khapugin2019local}. In this problem, a company enters the market (referred to as the follower in the CFLP context) with the objective of locating $r$ facilities to maximize its market share (capturing customer demand). Customers have associated demand for various commodities, which can be fulfilled from different facilities. Subsequently, customers determine their routing for collecting their demanded commodities from both existing and newly located facilities. In this problem, the entering firm acts as the leader, while the customers serve as the followers.\\

To solve the problem, we present a reformulation of the bilevel problem, establishing an equivalence with the maximal covering location problem. This reformulation enables exact and efficient problem-solving. Computational experimentation validates this approach.\\



%% file: abstracts/tex/IWOLOCA24_Abstract_7279.tex



\index[author]{Filippi, Carlo \\\hspace{0.5cm}
	 University of Brescia, Italy, \email{carlo.filippi@unibs.it}}
\index[author]{Guastaroba, Gianfranco \\\hspace{0.5cm}
	 University of Brescia, Italy, \email{gianfranco.guastaroba@unibs.it}}
\index[author]{Salazar-Gonz\'alez, Juan-Jos\'e \\\hspace{0.5cm}
	 University of La Laguna, Spain, \email{jjsalaza@ull.es}}


\title{Cross-docking platforms design and \newline mixed binary quadratic model \newline for distributionally robust optimization}


\author{
   Carlo Filippi\,\textsuperscript{a}, Gianfranco Guastaroba\,\textsuperscript{a}, Juan-José Salazar-González\,\textsuperscript{c,*}}

\address{
      \textsuperscript{a* } DEM, University of Brescia, Brescia, Italy, \email{carlo.filippi,gianfranco.guastaroba\}@unibs.it}\\
      \textsuperscript{b } IMAULL, University of La Laguna, La Laguna 38200, Tenerife, Spain \email{jjsalaza@ull.es}\\
      \textsuperscript{c } Applied Mathematics Dep., Universidad del Pa\'is Vasco
(UPV/EHU), Bilbao (Bizkaia), Spain\email{: aitziber.unzueta@ehu.eus}\\
}


\begin{abstract_online}
{Solving a multi-source capacitated facility location problem with a fairness objective}
{
 Carlo Filippi\,\textsuperscript{a}, Gianfranco Guastaroba\,\textsuperscript{a}, Juan-José Salazar-González\,\textsuperscript{b,*}
}
{
\textsuperscript{a }DEM, University of Brescia, Brescia, Italy, \\ 
\email{carlo.filippi,gianfranco.guastaroba\}@unibs.it}\\
\textsuperscript{b }IMAULL, University of La Laguna, La Laguna 38200, Tenerife, Spain \email{jjsalaza@ull.es}\\
\textsuperscript{*} Presenting author. \\
}
{
Location, Fairness, Cutting plane methods
}
\end{abstract_online}

	\section*{Introduction}
	In several application domains, a decision-maker has to locate a set of facilities, which provide a service (or good), and that are reached by ``customers'', usually individuals, at their own costs. In these situations, improving \emph{system efficiency} requires minimizing the average cost paid by customers to reach a facility, while improving \emph{fairness} requires minimizing the variability in the cost distribution  \cite{2021-FilGuaSpe}. Further, efficiency and fairness are commonly recognized as two competing objectives, where the former is usually the primary objective in the private sector, whereas the latter is normally the primary objective in the public sector \cite{eiselt1995objectives}.
	
	In this paper, we measure \emph{accessibility fairness} by using the \textit{conditional $\beta$-mean} \cite{2002-OgrZaw}, an equity measure strictly related to the Conditional Value-at-Risk, a popular risk measure \cite{2019-FilOgrSpe}. 
	
	As long as single-source location problems are involved (i.e., problems where the demand from a customer cannot be split among different facilities), the conditional $\beta$-mean can be incorporated into a MILP model preserving linearity \cite{2021-FilGuaSpe}.
	On the contrary, in multi-source location problems where the demand of a given customer can be allocated to more than one facility, incorporating the conditional $\beta$-mean into a MILP may lead to mixed integer non-linear programs. An example can be found in \cite{2018-ChaMit}.
	
	Our contributions can be summarized as follows: (1) we provide a non-linear formulation for the fair multi-source capacitated facility location problem;
	(2) we propose two linear reformulations of the latter model by using ideas from bilevel programming;
	(3) we develop solution methods that work along the general lines of a cutting-plane algorithm;
	(4) we show the effectiveness of the proposed algorithms.
	
	\section*{Models}
	The multi-source capacitated facility location problem is formulated as follows:
	\begin{alignat}{3}
		(\text{\ttfamily{MSCFL}}) \quad \min \quad &  \sum_{j \in J} f_j y_j + \sum_{i \in I} \sum_{j \in J} c_{ij} d_ i x_{ij} \label{MSCFL-obj} \\
		\text{subject to} \quad & \sum_{j \in J} x_{ij} = 1 & (i \in I) \label{MSCFL-assignment}\\
		&  \sum_{i \in I} d_i x_{ij} \leq s_j y_j  & (j \in J) \label{MSCFL-fullcapacity}\\
		&  y_{j} \in \{0,1\} & (j \in J) \label{MSCFL-ybinary}\\
		&  x_{ij} \geq 0 & (i \in I;  j \in J).
		\label{MSCFL-x>0}
	\end{alignat}%
	where $I$ is the set of customers, $J$ is the set of potential facilities, $d_i$ is the demand of customer $i \in I$, $f_j$ and $s_j$ are, respectively, the fixed cost and the capacity of facility $j \in J$, and finally $c_{ij}$ is the distance between $i$ and $j$. Let $\mathcal{XY}$ denote the feasible set of \text{\ttfamily{MSCFL}}.
	
	\begin{Def}
		Given a parameter $\beta \in (0,1]$ and a vector $(\mathbf{x}, \mathbf{y}) \in \mathcal{XY}$, the \emph{conditional $\beta$-mean} $M_{\beta}(\mathbf{x})$ is the average distance travelled by the $100 \times \beta$ percent of total demand that travels the longest distance to reach the assigned facility.
		\label{def:CBM}
	\end{Def}
	
	The conditional $\beta$-mean for a \text{\ttfamily{MSCFL}} problem can be expressed as:
	{\footnotesize
		\begin{equation}\label{eq:dualCBM}
			M_{\beta}(\mathbf{x}) = \min\left\{D_{\beta}u + \sum_{i \in I} \sum_{j \in J} d_i v_{ij} x_{ij}  \bigg\rvert u + v_{ij} \geq \frac{c_{ij}}{D_{\beta}};\ v_{ij}\geq 0\ (i \in I,j \in J)\right\}
		\end{equation}
	}%
	where $D_{\beta} = \beta\cdot\sum_{i \in I}d_i$. Problem (\ref{eq:dualCBM}) can be embedded into \text{\ttfamily{MSCFL}}, obtaining a non-linear formulation for the \textit{fair} \text{\ttfamily{MSCFL}} \textit{problem}. 
	%
	We linearise such a  formulation by transforming it into an equivalent bilevel program, and then employing standard techniques to reduce the bilevel program into a single-level optimization program. In this way, we produce two different formulations and discuss computational results.

%% file: abstracts/tex/IWOLOCA24_Abstract_7300.tex


\index[author]{Gar\'in, Mar\'ia Araceli\\\hspace{0.5cm}
        Universidad del Pa\'is Vasco, Spain, \email{mariaaraceli.garin@ehu.eus}}
\index[author]{Escudero, Laureano F.\\\hspace{0.5cm}
          Universidad Rey Juan Carlos, Spain,  \email{laureano.escudero@urjc.es}}
\index[author]{Unzueta, Aitziber\\\hspace{0.5cm}
         Universidad del Pa\'is Vasco, Spain, \email{aitziber.unzueta@ehu.eus}}


\title{Cross-docking platforms design and \newline mixed binary quadratic model \newline for distributionally robust optimization}


\author{
   María Araceli Gar{\'i}n\,\textsuperscript{a,*}, Laureano F. Escudero \,\textsuperscript{b}, Aitziber Unzueta \,\textsuperscript{c}}

\address{
      \textsuperscript{a } Quantitative Methods Dep., Universidad del Pa\'is Vasco
(UPV/EHU), Bilbao (Bizkaia), Spain\email{:  mariaaraceli.garin@ehu.eus}\\
      \textsuperscript{b }Statistics and Operations Research Area, Universidad Rey Juan Carlos
(URJC), Móstoles (Madrid), Spain\email{: laureano.escudero@urjc.es}\\
      \textsuperscript{c } Applied Mathematics Dep., Universidad del Pa\'is Vasco
(UPV/EHU), Bilbao (Bizkaia), Spain\email{: aitziber.unzueta@ehu.eus}\\
}


\begin{abstract_online}
{Cross-docking platforms design and mixed binary quadratic model for distributionally robust optimization}
{
  María Araceli Gar{\'i}n\,\textsuperscript{a,*}, Laureano F. Escudero\,\textsuperscript{b}, Aitziber Unzueta\,\textsuperscript{c}
}
{
\textsuperscript{a }Quantitative Methods Dep., Universidad del Pa\'is Vasco (UPV/EHU), Bilbao (Bizkaia), Spain, \email{mariaaraceli.garin@ehu.eus}\\
\textsuperscript{b }Statistics and Operations Research Area, Universidad Rey Juan Carlos (URJC), Móstoles (Madrid), Spain, \email{laureano.escudero@urjc.es}\\
\textsuperscript{c }Applied Mathematics Dep., Universidad del Pa\'is Vasco (UPV/EHU), Bilbao (Bizkaia), Spain, \email{aitziber.unzueta@ehu.eus}\\
\textsuperscript{*} Presenting author. \\
}
{
Cross-dock door design, stochastic two-stage scenario setting, mixed binary quadratic, distributionally robust optimization
}
\end{abstract_online}


The Cross-dock Door Design Problem (CDDP) consists of deciding on the number and capacity of inbound and outbound doors,
to minimize the construction and exploitation cost of the
infrastructure. The uncertainty, realized in a two-stage scenario
setting, lies in the occurrence of these nodes, the delivering
material volume and cost, and the possible capacity’s disruption of the doors.
A mixed 0-1 quadratic model for distributionally
  robust optimization is introduced.

\section*{Introduction}
Distributionally robust optimization (DRO) is motivated as a counterpart of the usually unknown underlying probability distribution (PD) followed by the uncertainty in dynamic problems. An approach is presented for the highly combinatorial Cross-dock Door Design Problem (CDDP) solving to decide the number and nominal capacity of the strip and stack doors.
It is assumed that the uncertainty is represented in the second stage of a two-stage setting,
where a set of scenarios for the uncertain parameters is considered
for each member of an ambiguity set that has been previously
generated.

Initially, a mixed binary quadratic DRO risk-neutral (RN) model is presented.
The objective function to minimize is included by the the cross-dock
infrastructure building cost (i.e., first stage door selection and
installation cost) plus the highest second stage assignment expected
cost in the scenarios, among the ambiguity set members.
The constraint system is composed of the first stage constraints and
the
second stage constraint set for each member of the ambiguity set.

In a second step, a risk-averse (RA) formulation is developed
 by adding the system of constraints and variables related to the
 risk-averse measure called second-order stochastic dominance to
 prevent, up to some extent, the negative implication of black swan
 scenarios in the value of the objective function to minimize.

Given the difficulty on solving this highly combinatorial problem, a mathematically equivalent MILP formulation is considered.
Since, even this formulation is still difficult to solve, a  scenario
Cluster Lagrangean Decomposition (CLD) is introduced for obtaining lower bounds to the MILP model. A lazy scheme for obtaining feasible solutions to the original model based on the CLD one is considered, by exploiting the special structure of the problem. A computational study validates the proposal.

%% file: abstracts/tex/IWOLOCA24_Abstract_9592.tex


\index[author]{Landete Ruiz, Mercedes \\\hspace{0.5cm}
	 Universidad Miguel Hern\'andez de Elche, Spain, \email{landete@umh.es}}
\index[author]{Muñoz-Ocaña, Juan Manuel \\\hspace{0.5cm}
	 Universidad de C\'adiz, Spain, \email{juanmanuel.munoz@uca.es}}
\index[author]{Rodr\'iguez-Ch\'ia, Antonio M. \\\hspace{0.5cm}
	 Universidad de C\'adiz, Spain, \email{antonio.rodriguezchia@uca.es}}
\index[author]{Saldanha-da-Gama, Francisco \\\hspace{0.5cm}
	 Sheffield University, United Kingdom, \email{francisco.saldanha-da-gama@sheffield.ac.uk}}


\title{Formulations and resolution procedures for upgrading hub networks }


\author{
 Mercedes Landete Ruiz\,\textsuperscript{a}, Juan M. Muñoz-Ocaña\,\textsuperscript{b,*}, Antonio M. Rodr\'iguez-Ch\'ia\,\textsuperscript{b}, Francisco Saldanha-da-Gama\,\textsuperscript{c}}

\address{
      \textsuperscript{a }Universidad Miguel Hern\'andez,  \email{landete@umh.es}\\
      \textsuperscript{b }Universidad de C\'adiz, \email{juanmanuel.munoz@uca.es}, \email{antonio.rodriguezchia@uca.es}\\
      \textsuperscript{c }Sheffield University Management School, \email{francisco.saldanha-da-gama@sheffield.ac.uk}\\
}


\begin{abstract_online}
{Formulations and resolution procedures for upgrading hub networks}
{
 Mercedes Landete Ruiz\,\textsuperscript{a}, Juan M. Muñoz-Ocaña\,\textsuperscript{b,*}, Antonio M. Rodr\'iguez-Ch\'ia\,\textsuperscript{b}, Francisco Saldanha-da-Gama\,\textsuperscript{c}
}
{
      \textsuperscript{a }Universidad Miguel Hern\'andez,  \email{landete@umh.es}\\
      \textsuperscript{b }Universidad de C\'adiz, \email{juanmanuel.munoz@uca.es}, \email{antonio.rodriguezchia@uca.es}\\
      \textsuperscript{c }Sheffield University Management School, \email{francisco.saldanha-da-gama@sheffield.ac.uk}\\
      \textsuperscript{*} Presenting author. \\
}
{
Hub location, Edge upgrading
}
\end{abstract_online}


\section*{Introduction}
Hub location problems have become an important stream of research within Location Science due to their relevance in many applications emerging in logistics, telecommunications, and transportation \cite{Alumur}. This talk presents different formulations for the single-allocation hub location problem with edge upgrades. 
In the broad context of network design and optimization, upgrading aims to reduce the transportation costs between sites that are connected by the upgraded edges. Two types of connections are considered to be upgraded: inter-hubs edges and edges that connect origin/destination sites to hubs.

Upgrading may be motivated by different contexts.
In telecommunications, upgrading an edge may involve using a higher speed cable between a pair of servers.
Similarly, in logistics and transportation, upgrading an edge may correspond to using a faster transportation mode, such as using a plane instead of a truck between a pair of cities, thus achieving a reduced travel time.
Another possibility is to redesign a trajectory (e.g. using a highway instead of a secondary road).

\section*{Formulations for solving the problem}
We consider a maximum number of connections that can be upgraded between hubs and a maximum number of connections between spokes (non-hubs) and hubs.
This may correspond, for instance, to a limited number of alternative vehicles available for accomplishing the improvement.

We start by considering complete hub networks.
Afterward, we extend the work to problems in which hub network design decisions have also to be made. 
For the latter, no specific topology is imposed for the hub-level network.
Given that the objective function accounts only for the transportation costs, a constraint is set on the number of hubs to install \cite{ContrerasOKelly:2019}.
These formulations present high computing burden, and therefore, various valid inequalities are included as cuts in branch-and-cut procedures.
The methodological developments are computationally tested for two purposes: firstly, to emphasize the relevance of embedding upgrading decisions in single-allocation hub location problems by analyzing the savings achieved and the structural changes in the hub networks after upgrading edges; secondly, to evaluate the effectiveness of the formulations and solution methods for different instances.


%% file: abstracts/tex/IWOLOCA24_Abstract_5088.tex


\index[author]{L\'opez-S\'anchez, Aitor \\\hspace{0.5cm}
	 Universidad Rey Juan Carlos, Spain, \email{aitor.lopez@urjc.es}}
\index[author]{Lujak, Marin \\\hspace{0.5cm}
	 Universidad Rey Juan Carlos, Spain, \email{marin.lujak@urjc.es}}
\index[author]{Semet, Frederic \\\hspace{0.5cm}
	 Centrale Lille, France, \email{frederic.semet@centralelille.fr}}
\index[author]{Billhardt, Holger \\\hspace{0.5cm}
	 Universidad Rey Juan Carlos, Spain, \email{holger.billhardt@urjc.es}}


\title{Synchronization routing for agricultural \\ vehicles and implements}


\author{
Aitor L\'opez-S\'anchez\,\textsuperscript{a,b,*}, Marin Lujak\,\textsuperscript{a}, Fr\'ed\'eric Semet\,\textsuperscript{b}, Holger Billhardt\,\textsuperscript{a}}

\address{
      \textsuperscript{a }CETINIA, University Rey Juan Carlos, Madrid Spain,  \email{aitor.lopez@urjc.es}, \\  \email{marin.lujak@urjc.es},  \email{holger.billhardt@urjc.es}\\
      \textsuperscript{b }UMR 9189 - CRIStAL, Centrale Lille, Univ. Lille, CNRS, Inria, France, \\ \email{frederic.semet@centralelille.fr}\\
}


\begin{abstract_online}
{Synchronization routing for agricultural vehicles and implements}
{
Aitor L\'opez-S\'anchez\,\textsuperscript{a,b,*}, Marin Lujak\,\textsuperscript{a}, Fr\'ed\'eric Semet\,\textsuperscript{b}, Holger Billhardt\,\textsuperscript{a}
}
{
\textsuperscript{a }CETINIA, University Rey Juan Carlos, Madrid, Spain,  \email{aitor.lopez@urjc.es}, \\  \email{marin.lujak@urjc.es},  \email{holger.billhardt@urjc.es}\\
\textsuperscript{b }UMR 9189 - CRIStAL, Centrale Lille, Univ. Lille, CNRS, Inria, France, \\ \email{frederic.semet@centralelille.fr}\\
\textsuperscript{*} Presenting author. \\
}
{
Routing, Synchronization, Column generation, Agriculture
}
\end{abstract_online}


\section*{Introduction}

Agriculture technology is experiencing a revolution to enhance efficiency, sustainability, and productivity. Autonomous agriculture mobile robots (agribot) are increasingly used.  Current approaches for Agricultural Vehicle Routing Problems (AVRPs) focus on homogeneous fleets, dealing with a single type of task and crop \cite{utamima2022agricultural}. In contrast, real-world farming requires coordination between different tractors/agribots and implements to manage different crops and tasks, such as plowing, fertilizing, and harvesting. The Agricultural Fleet Vehicle Routing Problem with Implements (AFVRPI) coordinates the routes for mixed fleets, optimizing the cost for covering task execution and synchronization in modern agriculture.

Movement synchronization between two vehicle classes (tractors and implements) is necessary, where the tractors have autonomous movement capacity and implements depend on the tractors to move \cite{soares2023synchronisation}.
Various approaches in the literature address similar challenges. The Vehicle Routing Problem with Trailers and Transshipments \cite{drexl2013applications}  allows for the detachment and reattachment of trailers between trucks. The Active Passive Vehicle Routing Problem introduces a scenario with active and passive vehicles, where the active vehicles displace the passive ones \cite{tilk2018branch}. Lastly, one of the first movement synchronization constraints in precision agriculture is the Synchronized Sprayer Tanker Routing Problem with Variable Service Time, which involves coordination between tender tankers and sprayer fleets \cite{alkaabneh2024matheuristic}.

\section*{Problem definition and solution approach}

The Agricultural Fleet Vehicle Routing Problem with Implements (AFVRPI) is composed of a fleet comprising both tractors and implements $(\mathcal{F} = \mathcal{V} \cup \mathcal{M})$, and their compatibilities with each other. The transportation network is modeled by a directed graph $(\mathcal{N}, \mathcal{A})$. Nodes consist of four distinct sets: $\mathcal{N}_{tasks}$ agricultural tasks, $\mathcal{N}_{depots}$ tractor and implement depots, $\mathcal{N}_{detach}$ detaching nodes and $\mathcal{N}_{attach}$ attaching nodes. Each task has a given demand, service time, and time window.
Arcs denote spatial and temporal connectivity and let us denote transfer arcs as $(d, a) \in \mathcal{A}_{transfer}$, including arcs where implements can be detached $d \in \mathcal{N}_{detach}$/attached $a \in \mathcal{N}_{attach}$ to tractors.  The goal of the AFVRPI is to find a set of compatible routes for tractors and implements that visit all the tasks minimizing the overall cost and respecting the movement constraints.

The proposed solution approach is composed of an extended master problem formulation with independent subproblems associated with each tractor and implement, and solved with column generation, which enables distributed and asynchronous problem-solving.
Tractor subproblems are the  Elementary Shortest Path Problem with Resource Constraints (ESPPRC) incorporating linear costs \cite{tilk2018branch} with two resources: the distance, restricted by vehicle autonomy, and the time, constrained by task time windows. Implement subproblems are also the ESPPRC considering only demand constraints, limited by its capacities.  
Preliminary computational results demonstrate the efficiency of this approach.


%% file: abstracts/tex/IWOLOCA24_Abstract_8729.tex


\index[author]{Melero, Inigo Martin \\\hspace{0.5cm}
	 Universidad Miguel Hernandez de Elche, Spain, \email{inigo.martin@goumh.umh.es}}
\index[author]{Landete Ruiz, Mercedes \\\hspace{0.5cm}
	 Universidad Miguel Hern\'andez de Elche, Spain, \email{landete@umh.es}}
\index[author]{Alcaraz Soria, Javier \\\hspace{0.5cm}
	 Universidad Miguel Hernandez de Elche, Spain, \email{jalcaraz@umh.es}}


\title{The Effect of Budget Limiting on the Linear Ordering Problem}


\author{
 Inigo Martin Melero\,\textsuperscript{a,b,*}, Mercedes Landete Ruiz\,\textsuperscript{b}, Javier Alcaraz Soria\,\textsuperscript{b}}

\address{
      \textsuperscript{a }European Organisation for Nuclear Research (CERN),  \email{inigo.martin.melero@cern.ch}\\
      \textsuperscript{b }Universidad Miguel Hernandez de Elche, \email{inigo.martin@goumh.umh.es}, \email{landete@umh.es},\\
      \email{jalcaraz@umh.es}\\
}


\begin{abstract_online}
{The Effect of Budget Limiting on the Linear Ordering Problem}
{
 Inigo Martin Melero\,\textsuperscript{a,b,*}, Mercedes Landete Ruiz\,\textsuperscript{b}, Javier Alcaraz Soria\,\textsuperscript{b}
}
{
\textsuperscript{a }European Organisation for Nuclear Research (CERN),  \email{inigo.martin.melero@cern.ch}\\
\textsuperscript{b }Universidad Miguel Hernandez de Elche, \email{inigo.martin@goumh.umh.es}, \email{landete@umh.es},\\
\email{jalcaraz@umh.es}\\
\textsuperscript{*} Presenting author. \\
}
{
Linear Ordering Problem, Bilevel Programming, Routing
}
\end{abstract_online}


\section*{Introduction}

Routing, inside the category of location problems, can be interpreted and modelled as ranking problems: the first locations to be visited are those with the best ranking scores and the direction of the route indicates the position in the ranking. Therefore, the position that an element occupies in a route or ranking usually depends on different criteria. 

Considering a given route and focusing on one of the visited elements in the route that is not in first place and that has budget to improve its position, in this work we model the problem of deciding in which criteria it should improve in order to be visited as soon as possible. We assume that the order of the route is calculated by solving the Linear Ordering Problem \cite{Alcaraz} that makes use of the information in the different criteria. We propose a bilevel model in which the first level improves the position of the selected element and the second level solves the Linear Ordering Problem. We propose a genetic heuristic resolution algorithm and compare the results with the exact resolution for instances where the exact can be solved. We analyze the key properties of the problem and the optimization model.

%% file: abstracts/tex/IWOLOCA24_Abstract_2730.tex


\index[author]{Piedra-de-la-Cuadra, Ram\'on \\\hspace{0.5cm}
	 Universidad de Sevilla, Spain, \email{rpiedra@us.es}}
\index[author]{Ortega, Francisco A. \\\hspace{0.5cm}
	 Universidad de Sevilla, Spain, \email{riejos@us.es}}


\title{Optimization approaches to scheduling working shifts for train dispatchers}


\author{
 Ram\'on Piedra-de-la-Cuadra\,\textsuperscript{a,*} , Francisco A. Ortega\,\textsuperscript{b}}

\address{
      \textsuperscript{a }Universidad de Sevilla,  \email{rpiedra@us.es}\\
      \textsuperscript{b }Universidad de Sevilla, \email{riejos@us.es}\\
}


\begin{abstract_online}
{Optimization approaches to scheduling working shifts for train dispatchers}
{
 Ram\'on Piedra-de-la-Cuadra\,\textsuperscript{a,*}, Francisco A. Ortega\,\textsuperscript{a}
}
{
\textsuperscript{a }Universidad de Sevilla,  \email{rpiedra@us.es}, \email{riejos@us.es}\\
\textsuperscript{*} Presenting author. \\
}
{
Railway dispatching, Shift scheduling, Integer Programming
}
\end{abstract_online}


\section*{Introduction}

Scheduling work shifts involves the strategic assignment of tasks, the determination of their locations, start and end times and, accordingly, the allocation of personnel. For railway dispatchers, this responsibility is paramount as they plan, coordinate, and oversee train movements within the rail network, ensuring safety and efficiency. The complexity of this task is magnified by various constraints such as legal limitations on shift duration and on the consecutive night shifts, in addition to other operational restrictions like the number of areas a dispatcher can cover, and the need to avoid undesirable schedules.
Traditionally, train dispatcher schedules have been manually crafted, which is time-consuming, and it can be prone to errors given the high number of constraints to be considered. To address this, an optimization model is proposed to automate the scheduling process. 
By automating shift scheduling using optimization models, railway companies can streamline operations, reduce administrative burden, and improve overall efficiency. Moreover, by considering both legal requirements and employee preferences, these automated systems can foster a better work-life balance and enhance employee morale. In conclusion, the implementation of optimization models offers a promising approach to address the complexities of scheduling work shifts for train dispatchers while ensuring safety, efficiency, and the employee well-being.\\
Our proposed model considers various legal limitations, assignments of feasible areas, and the quality of shifts to avoid undesirable configurations, even if permissible by law or union agreements.
Input data for the model includes geographical areas within the railway network, feasible combinations of zones based on the existing infrastructure, and initial assignments for dispatchers according to their preferences. Additionally, constraints such as task load per dispatcher, shift duration limits, rest periods, and controllable areas per dispatcher are incorporated. The objective is to generate schedules that minimize the number of dispatchers while adhering to legal and operational constraints and optimizing employee satisfaction.


%% file: abstracts/tex/IWOLOCA24_Abstract_1386.tex


\index[author]{Puerto, Justo \\\hspace{0.5cm}
	 Universidad de Sevilla, Spain, \email{puerto@us.es}}
\index[author]{Valverde, Carlos \\\hspace{0.5cm}
	 Universidad de Sevilla, Spain, \email{cvalverde@us.es}}


\title{The Hampered K-Median Problem with Neighbourhoods}


\author{
 Justo Puerto \textsuperscript{a} and Carlos Valverde\textsuperscript{a,*} }

\address{
      \textsuperscript{a }Instituto de Matemáticas, Universidad de Sevilla,  \email{puerto@us.es},   \email{cvalverde@us.es}\\
      \textsuperscript{*} Presenting author. \\
}

\keywords{Facility location, Continuous location, Barriers, Mixed integer Conic programming}

\begin{abstract_online}
{The Hampered K-Median Problem with Neighbourhoods}
{
Justo Puerto\,\textsuperscript{a}, Carlos Valverde\,\textsuperscript{a,*}
}
{
\textsuperscript{a }Instituto de Matemáticas, Universidad de Sevilla,  \email{puerto@us.es},   \email{cvalverde@us.es}\\
\textsuperscript{*} Presenting author. \\
}
{
Facility location, Continuous location, Barriers, Mixed integer Conic programming
}
\end{abstract_online}


\section*{Abstract}
This paper deals with facility location problems in a continuous space with neighbours and barriers. Each one of these two elements, neighbours and barriers, makes the problems harder than their standard counterparts. Combining all together results in a new challenging problem that, as far as we know, has not been addressed before, but has applications for inspection and surveillance activities and the delivery industry assuming uniformly distributed demand in some regions. Specifically, we analyse the $k$-Median problem with neighbours and polygonal barriers in two different situations. None of these problems can be  seen as a simple incremental contribution since in both cases the tools required to analyse and solve them go beyond any standard overlapping of techniques used in the separated problems. As a first building block, we deal with the problem assuming that the neighbourhoods are not visible from one another and therefore there are no rectilinear paths that join two of them without crossing barriers. Under this hypothesis, we derive a valid mixed-integer linear formulation. Removing that hypothesis leads to a more general and realistic problem, but at the price of making it more challenging. Adapting the elements of the first formulation, we also develop another valid mixed-integer bilinear formulation.  Both formulations rely on tools borrowed from computational geometry that allow to handle polygonal barriers and neighbours that are second-order cone (SOC) representable, which we preprocess and strengthen with valid inequalities. These mathematical programming formulations are also instrumental to derive an adapted matheuristic algorithm that provides good quality solutions for both problems in short computing time. The paper also reports extensive computational experience, counting 2400 experiments, showing that our exact and heuristic approaches are useful: the exact approach can solve optimally instances with up to 50 neighbourhoods and different number of barriers within one hour of CPU time, whereas the matheuristic approach always returns good feasible solutions in less than 300 seconds.

%% file: abstracts/tex/IWOLOCA24_Abstract_3749.tex


\index[author]{Segura, Paula \\\hspace{0.5cm}
	 Universitat Politècnica de València, Spain, \email{psegmar@upvnet.upv.es}}
\index[author]{Plana, Isaac \\\hspace{0.5cm}
	 Universitat de València, Spain, \email{isaac.plana@uv.es}}
\index[author]{Sanchis, José María\\\hspace{0.5cm}
	 Universitat Politècnica de València, Spain, \email{jmsanchis@mat.upv.es}}


\title{Exact approaches for the Chinese Postman Problem with load–dependent costs}


\author{
 Paula Segura\,\textsuperscript{a,*} , Isaac Plana\,\textsuperscript{b}, José María Sanchis\,\textsuperscript{a}}

\address{
\textsuperscript{a }Universitat Politècnica de València,  \email{psegmar@upvnet.upv.es},   \email{jmsanchis@mat.upv.es}\\
\textsuperscript{b }Universitat de València, \email{isaac.plana@uv.es}\\
\textsuperscript{*} Presenting author. \\
}


\begin{abstract_online}
{Exact approaches for the Chinese Postman Problem with \newline load–dependent costs}
{
Paula Segura\,\textsuperscript{a,*}, Isaac Plana\,\textsuperscript{b}, José María Sanchis\,\textsuperscript{a}
}
{
\textsuperscript{a }Universitat Politècnica de València,  \email{psegmar@upvnet.upv.es},   \email{jmsanchis@mat.upv.es}\\
\textsuperscript{b }Universitat de València, \email{isaac.plana@uv.es}\\
\textsuperscript{*} Presenting author. \\
}
{
arc routing problems, load-dependent costs, mathematical formulations, branch and cut
}
\end{abstract_online}


\section*{Introduction}

The Chinese Postman Problem (CPP) is a classical arc routing problem whose goal is to find a minimum-cost tour on a connected undirected graph that traverses each edge at least once (\cite{Laporte14}). In transportation systems, the level of emissions from a vehicle is influenced by factors beyond the distance traveled, such as its load. Motivated from the desire to reduce pollution, the Chinese Postman Problem with load–dependent costs (CPP-LC) was introduced in \cite{CELPS} as an extension of the CPP in which the cost of traversing an edge depends on its length and also on the weight of the vehicle at the moment the edge is traversed. In this talk, we summarize the two mathematical programming formulations proposed in the literature for the CPP-LC and present a new mixed-integer linear programming formulation for the problem, proposing several families of valid inequalities to reinforce such a formulation. We design a new branch-and-cut algorithm for the CPP-LC solution based on this formulation that incorporates the separation of the valid inequalities proposed. Some computational results obtained with our new exact procedure are compared with those already existing in the literature.



%% file: abstracts/tex/IWOLOCA24_Abstract_7463.tex


\index[author]{Benati, Stefano\\\hspace{0.5cm}
	 Università degli Studi di Trento, Italy, \email{stefano.benati@unitn.it}}
\index[author]{Puerto, Justo \\\hspace{0.5cm}
	 Universidad de Sevilla, Spain, \email{puerto@us.es}}
\index[author]{Temprano, Francisco\\\hspace{0.5cm}
	 Universidad de Sevilla, Spain, \email{ftgarcia@us.es}}


\title{New advances in hypergraph structure analysis}


\author{
Francisco Temprano\,\textsuperscript{a,*} , Stefano Benati\,\textsuperscript{b}, Justo Puerto\,\textsuperscript{a}}

\address{
      \textsuperscript{a }Universidad de Sevilla,  \email{ftgarcia@us.es},   \email{puerto@us.es}\\
      \textsuperscript{b }Università degli Studi di Trento, \email{stefano.benati@unitn.it}\\
}


\begin{abstract_online}
{New advances in hypergraph structure analysis}
{
 Francisco Temprano\,\textsuperscript{a,*}, Stefano Benati\,\textsuperscript{b}, Justo Puerto\,\textsuperscript{a}
}
{
\textsuperscript{a }Universidad de Sevilla,  \email{ftgarcia@us.es},   \email{puerto@us.es}\\
      \textsuperscript{b }Università degli Studi di Trento, \email{stefano.benati@unitn.it}\\
      \textsuperscript{*} Presenting author. \\
}
{
Network location, Mathematical programming, Branch-and-Price
}
\end{abstract_online}


\section*{Abstract}

This talk deals with the problem of locating subsets of nodes highly con
nected over hypergraphs. Due to the lack of universal consensus and de
veloped theory in literature, we must come up with a formal definition aim
of the problem, in order to propose new optimization models to solve it.
A hypergraph is known to be a generalization of a graph that allows us
to represent a huge amount of real-life interactions between elements of a
data sample that the classic and original networks can not. Once we for
mally define the problem, we develop a large list of functions able to mea
sure the goodness of the node set subdivision. Next, we compare all these
methods by using of mathematical programming and extensive computa
tional experiments. Thus, we can conclude which methods describe better
the properties of partition structures and which ones perform better, since
multiple compact formulations and column generation algorithms are de
veloped to solve the partitioning hypergraph problem. Finally, we have
applied all methods to Eurobarometer data, showing the applicability of
the introduced methodology.

%% file: abstracts/tex/IWOLOCA24_Abstract_4806.tex


\index[author]{Unzueta, Aitziber\\\hspace{0.5cm}
         Universidad del Pa\'is Vasco, Spain, \email{aitziber.unzueta@ehu.eus}}
\index[author]{Escudero, Laureano F.\\\hspace{0.5cm}
          Universidad Rey Juan Carlos, Spain,  \email{laureano.escudero@urjc.es}}
\index[author]{Gar\'in, Mar\'ia Araceli\\\hspace{0.5cm}
        Universidad del Pa\'is Vasco, Spain, \email{mariaaraceli.garin@ehu.eus}}


\title{Cross-docking platforms design and \newline management under uncertainty}


\author{
 Aitziber Unzueta\,\textsuperscript{a,*} , Laureano F. Escudero\,\textsuperscript{b}, María Araceli Gar{\'i}n\,\textsuperscript{c}}

\address{
\textsuperscript{a* }Applied Mathematics Department, University of the Basque Country (UPV/EHU), Bilbao (Bizkaia), Spain\email{: aitziber.unzueta@ehu.eus}\\
\textsuperscript{b } Statistics and Operations Research Area, University King Juan Carlos (URJC), Móstoles (Madrid), Spain\email{: laureano.escudero@urjc.es}\\
\textsuperscript{c } Quantitative Methods Department, University of the Basque Country (UPV/EHU), Bilbao (Bizkaia), Spain\email{:  mariaaraceli.garin@ehu.eus}\\
\textsuperscript{*} Presenting author. \\
}


\begin{abstract_online}
{Cross-docking platforms design and management under uncertainty}
{
Aitziber Unzueta\,\textsuperscript{a,*}, Laureano F. Escudero\,\textsuperscript{b}, María Araceli Gar{\'i}n\,\textsuperscript{c}
}
{
\textsuperscript{a }Applied Mathematics Department, University of the Basque Country (UPV/EHU), Bilbao (Bizkaia), Spain, \email{aitziber.unzueta@ehu.eus}\\
\textsuperscript{b }Statistics and Operations Research Area, University King Juan Carlos (URJC), Móstoles (Madrid), Spain, \email{laureano.escudero@urjc.es}\\
\textsuperscript{c }Quantitative Methods Department, University of the Basque Country (UPV/EHU), Bilbao (Bizkaia), Spain, \email{mariaaraceli.garin@ehu.eus}\\
\textsuperscript{*} Presenting author. \\
}
{
Cross-dock door design, Two-stage stochastic quadratic
combinatorial optimization, Linearized mixed-integer programming, Constructive matheuristic.
}
\end{abstract_online}


The Cross-dock Door Design Problem (CDDP) consists of deciding on the number and capacity of inbound and outbound doors,
minimizing the construction and exploitation cost of the infrastructure. The uncertainty, realized in scenarios, lies in the occurrence of these nodes, the number and cost of the pallets, and the capacity’s disruption of the doors.
The CDDP is represented using a stochastic two-stage binary quadratic model (BQM).
Given the difficulty of solving this combinatorial problem, a mathematically equivalent MILP
formulation is introduced. Searching an optimal solution of the whole MILP problem is still impractical, thus,  a scenario cluster decomposition-based matheuristic algorithm is presented and a broad computational study to validate the proposal is carried out.

\section*{Introduction}

Given a network with a set of supplying (i.e. origin) nodes for different product types
and a set of receiving (i.e. destination) nodes for these products, usually a cross-dock entity may serve as a consolidation point. The origin nodes
can deliver the material at the cross-dock so that, after being classified by type and
destination, it can be transported to the destination nodes. A cross-dock infrastructure
has a number of inbound doors, which are called strip ones, and a number of outbound
doors, which are called stack ones, and each of them has a capacity for pallet handling
during a given time period.

Two main types of optimization problems arise in cross-dock dealing. One of them is the deterministic Cross-dock
Door Assignment Problem (CDAP), see \cite{Guignard20} and \cite{cdap23} for further information about solving methodologies. The other one is the stochastic  CDDP, the subject of this work.
Comprehensive reviews on uncertainty have been recently published in \cite{Ardakani20}.

In this work we introduce a two-stage stochastic BQP model for cross-dock door design planning as an extension of the deterministic CDAP to deal with the uncertainty. The first stage is devoted to the strategic decisions (i.e. the number of strip and stack doors and related nominal capacities), before the uncertainties related to the set of origin and destination nodes and the disruption of the doors capacity are unveiled. The second stage is devoted to the operational decisions (i.e. the assignment of doors to origin and destination nodes in the scenarios). The objective function to minimize is composed of the binary linear cross-dock door infrastructure construction cost and the binary quadratic cost function related to the CDAP scenario node-to-door assignment. We use the Reformulation Linearization Technique to develop a Linearized mixed Integer Programming problem (LIP) as a equivalent reformulation of the quadratic one.
We introduce two options for a scenario Cluster Decomposition (CD) of model LIP, where the special structure of the problem is benefitted from. The first option decomposes model LIP into scenario clusters for the first stage variables. The second option additionally does this for the strip and stack door problems.
Finally, we develop a matheuristic algorithm, based on feasible first stage solutions of the CD model, to obtain feasible solutions for the stochastic model.  It considers a linear search approach for large-scale instances to exploit the scenario structure of the second stage submodels. It has been computationally proved that the proposed scheme provides the same or a better incumbent solution value than the plain use of CPLEX and it requires a much smaller wall time.

%% file: 06_useful_info.tex
\section{Venues}

\textbf{Talks} will be held at two buildings: the \textbf{Instituto de Matemáticas de Granada (IMAG)} and \textbf{El Carmen de la Victoria}. 
\begin{itemize}
    \item The address of the \textbf{IMAG} is at Ventanilla nº 11 de Granada. The way of arriving at the institute is by Rector López-Argüeta Street. Follow the red path in the picture below.
    \begin{figure}[h]
	\centering	
		\fbox{\includegraphics[scale=0.4]{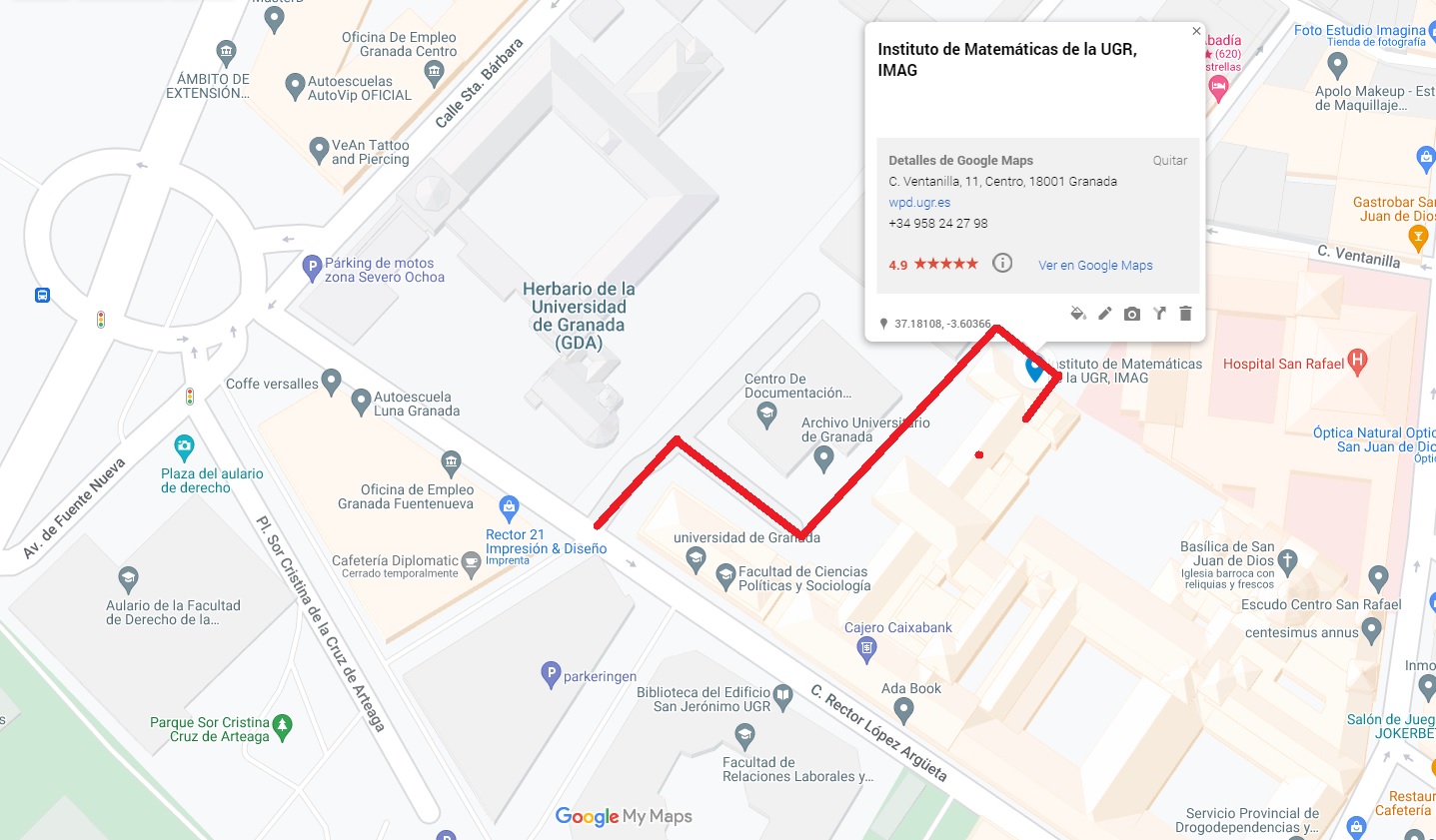}}
		\caption*{Access to IMAG.}
\end{figure}
    \item The address of \textbf{El Carmen de la Victoria} is at Encrucijada, 13. It is 30 minutes far away from the IMAG.

    Press the following link for directions in \textit{Google maps} from IMAG to El Carmen de la Victoria: \url{https://maps.app.goo.gl/ufN4UTY9kMEmQmKZA}
\end{itemize}

Our venues are at 1.8 km from each other. Commuting between them will take 28 min walking distance. Uber/Taxi is available in Granada, as well as urban bus (take one bus from Triunfo to Plaza Nueva and another from Plaza Nueva to Paseo de los Tristes).

\begin{center}
\fbox{\includegraphics[width=0.95\textwidth]{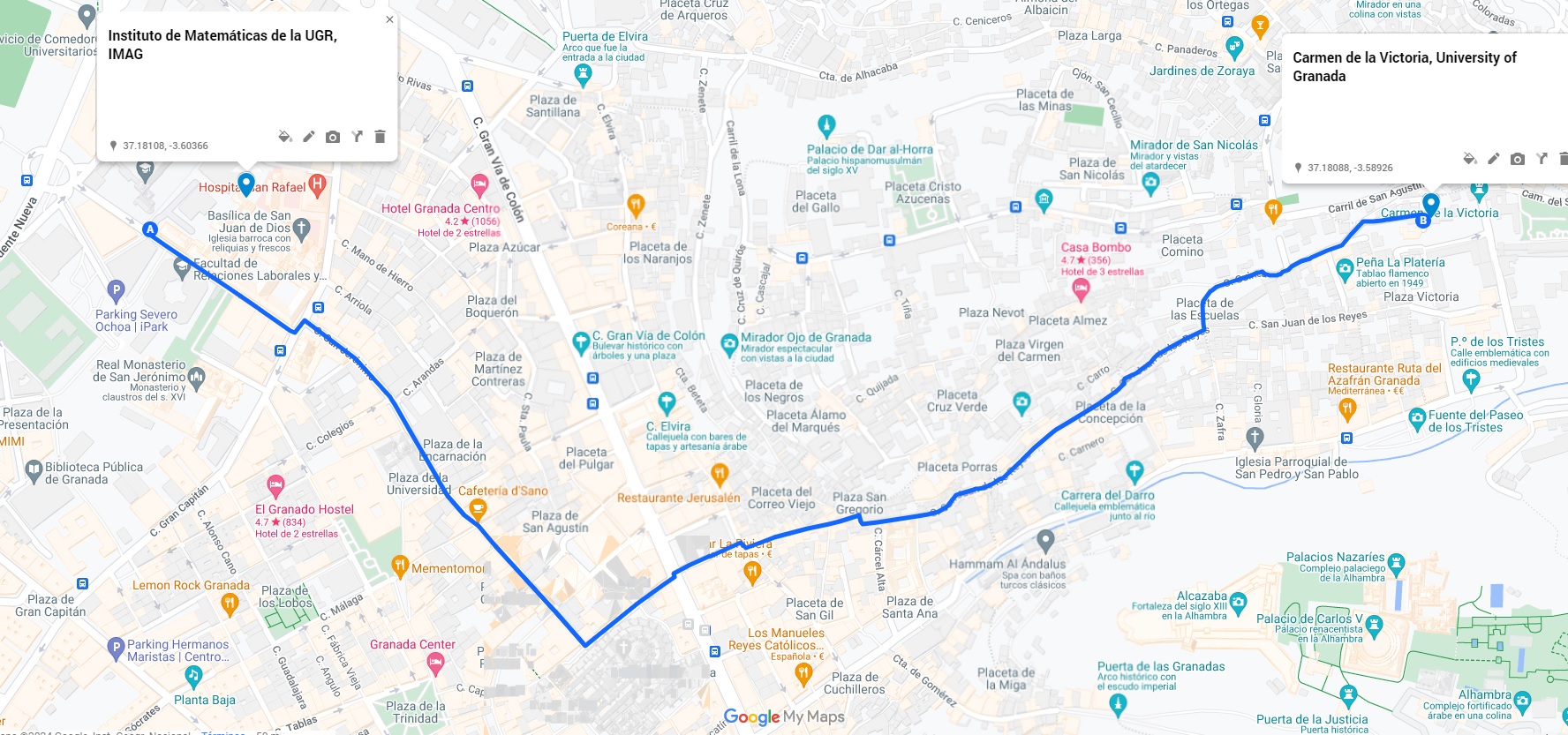}}
\end{center}

\section*{Wi-Fi}

Wi-Fi will be available during the conference using eduroam network.

\section{Social Activities}

\subsection*{$\star$ Wednesday, September 4, 2024: Welcome Reception}

On Wednesday, it will take place the welcome reception at 9 pm at {\color{\IScolor}\textbf{BHeaven}},  Calle Acera del Darro, 62, at the rooftop of Hotel Barceló Carmen Granada.

\subsection*{$\star$ Thursday, September 5, 2024: Guided Tour}

On Thursday, we will have a walk through Albayzin, the picturesque neighborhood, focusing on the viewpoints (\emph{miradores}). The tour will finish with a gastronomic tour to taste the specialty of Granada: the \emph{tapas}. The departure for this event will be at 6 pm, from El Carmen de la Victoria.

\subsection*{$\star$ Friday, September 6, 2024: Gala Dinner}

On Friday, our \textbf{Gala dinner} will be held at 9 pm at the gardens of {\color{\IScolor}\textbf{El Carmen de la Victoria}}.

\section*{Interactive map}

All locations related to the IWOLOCA 2024 conference can be found on this map, as well as points of interest in the city of Granada: 

\url{https://bit.ly/IWOLOCA24_spots}

\begin{center}
    \includegraphics[width=0.5\textwidth]{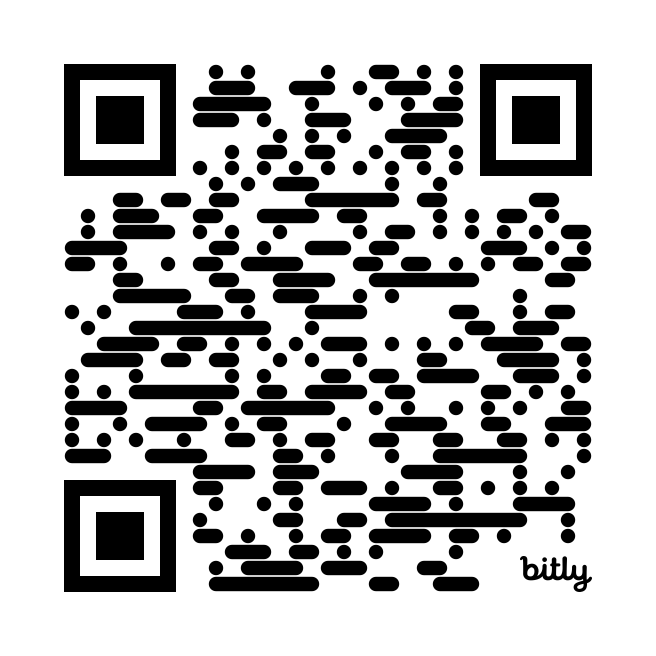}
\end{center}

%% file: 07_logos.tex
\begin{center}
The IWOLOCA is very grateful to the following sponsors.
\end{center}

\vfill

\begin{figure}[h]
	\includegraphics[width=\textwidth]{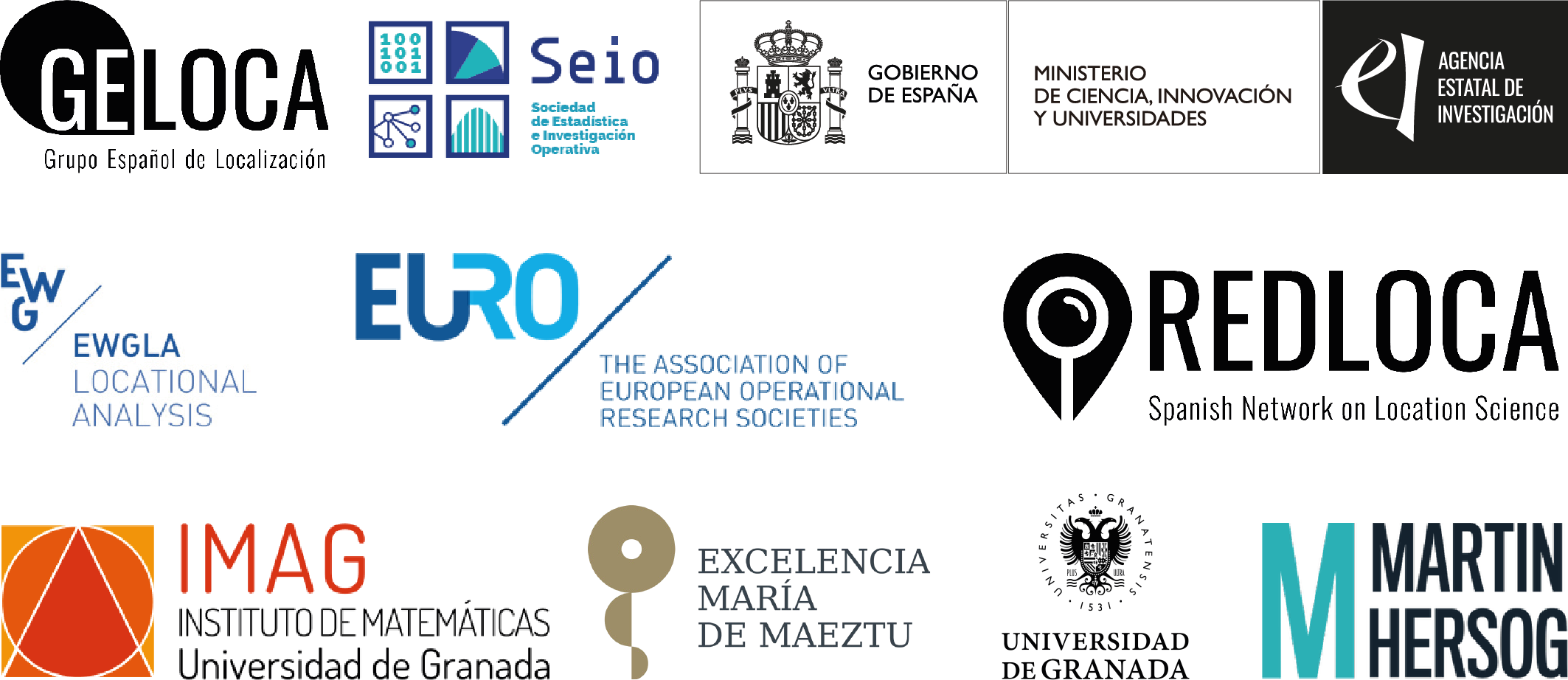}
\end{figure}
\vfill
XIII IWOLOCA 2024 has been partially supported by Grant RED2022-134149-T funded by MICIU/AEI /10.13039/501100011033 and IMAG-Maria de Maeztu grant CEX2020-001105-M /AEI /10.13039/501100011033